\documentclass[final]{siamart1116}
\usepackage[utf8]{inputenc}
\usepackage[english]{babel}
\usepackage[utf8]{inputenc} % Consente l'uso caratteri accentati italiani
\usepackage[a4paper,top=4cm,bottom=4cm,left=3cm,right=3cm]{geometry}

\usepackage{graphicx} %include i grafici
\usepackage{amssymb}
\usepackage{amsmath}

\usepackage{mathrsfs}
\usepackage{latexsym}
\DeclareGraphicsExtensions{.pdf,.jpg}
\usepackage{subfigure}
\usepackage{amsfonts}
\usepackage{bm}
\usepackage{verbatim}

\def \bmu {\bm{\mu}}
\def \uFE {u^\mathcal{N}(\bmu)}
\def \uRB {u^\mathcal{N}_N(\bmu)}
\def \pec {\mathbb{P}\mathtt{e}_K(x)}
\def \bbe {\bm{\beta}}
\def \R {\mathbb{R}}
\def \PGP {\textbf{PG }}
\def \PFS {\textbf{PFS }}

\newsiamremark{remark}{Remark}

\usepackage{lipsum}
\usepackage{amsfonts}
\usepackage{graphicx}
\usepackage{epstopdf}
\usepackage{algorithmic}
\ifpdf
  \DeclareGraphicsExtensions{.eps,.pdf,.png,.jpg, .JPG}
\else
  \DeclareGraphicsExtensions{.eps}
\fi

\numberwithin{theorem}{section}

\newcommand{\TheTitle}{\uppercase{Stabilized weighted reduced basis methods for parametrized advection dominated problems with random inputs} }
\newcommand{\ShortTitle}{Stabilized weighted RBM for random inputs advection dominated problems} 
\newcommand{\TheAuthors}{D. Torlo, F. Ballarin, and G. Rozza}

\headers{\ShortTitle}{\TheAuthors}

\title{{\TheTitle}\footnote[3]{ R\lowercase{eceived by the editors }J\lowercase{anuary 2, 2018; accepted for publication (in revised form) }S\lowercase{eptember 5, 2018; published
electronically }O\lowercase{ctober 25, 2018.}
\funding{T\lowercase{his work was funded by }E\lowercase{uropean }U\lowercase{nion }F\lowercase{unding for }R\lowercase{esearch and }I\lowercase{nnovation (project }H2020 ERC C\lowercase{o}G 2015 AROMA-CFD \lowercase{project 681447) and by the} INDAM-GNCS \lowercase{project.}}}}

\author{
  Davide Torlo \footnote[1]{\lowercase{math}L\lowercase{ab, }M\lowercase{athematics }A\lowercase{rea, }SISSA, V\lowercase{ia} B\lowercase{onomea 265,} I-34136 T\lowercase{rieste, }I\lowercase{taly}} \footnote[2]{ C\lowercase{urrent address:} I\lowercase{nstitut f\"ur }M\lowercase{athematik.} UZH, U\lowercase{niversit\"at }Z\lowercase{\"urich, }W\lowercase{interthurerstrasse 190,} CH-8057, Z\lowercase{\"urich,} S\lowercase{witzerland}}
  \and
  Francesco Ballarin\footnotemark[1]
   \and
 Gianluigi Rozza\footnotemark[1]
}

\usepackage{amsopn}

\begin{document}

\maketitle

\begin{abstract} 
In this work, we propose viable and efficient strategies for stabilized parametrized advection dominated problems, with random inputs. In particular, we investigate the combination of wRB (weighted reduced basis) method for stochastic parametrized problems with stabilized reduced basis method, which is the integration of classical stabilization methods (SUPG, in our case) in the Offline--Online structure of the RB method. Moreover, we introduce a reduction method that selectively enables online stabilization; this leads to a sensible reduction of computational costs, while keeping a very good accuracy with respect to high fidelity solutions. We present numerical test cases to assess the performance of the proposed methods in steady and unsteady problems related to heat transfer phenomena.
\end{abstract}

% REQUIRED

% REQUIRED

\section{Introduction}
Advection--diffusion equations are very important in many engineering applications, because they are used to model, for example, heat transfer phenomena \cite{incropera} or the diffusion phenomena, such as of pollutants in the atmosphere \cite{dede}. We are interested in studying related advection--diffusion PDEs when their P\'eclet numbers, representing, roughly, the ratio between the advection and the diffusion field, are high. Moreover, in such applications, we often need very fast evaluations of the approximated solution, depending on some input parameters, which may be deterministic or uncertain. This happens, for example, in the case of \emph{real-time} simulation or if we need to perform repeated approximations of solutions, for different input parameters. We find such \emph{many-query} situations in optimization problems, in which the objective function to be optimized depends on the parameters through the solution of a PDE or a system of PDEs.
\par The aim of this work is to study a stabilized reduced basis method suitable for the approximation of parametrized advection--diffusion partial differential equations (PDEs), in advection dominated cases, including a stochastic context, by considering random inputs. Indeed, the reduced basis (RB) method \cite{reduced_basis} has been devised to reduce the computational effort required by the repeated solution of parametrized problems. It provides rapidly approximation of solution of PDEs and it is able to guarantee the \emph{reliability} of the solution with a sharp and accurate \emph{a posteriori }error bound. In literature we can find many works about the application of the RB method to advection-diffusion problems, in particular with low P\'eclet number \cite{gelsomino, manzoni_45, rozza_nguyen}.
\par In contrast, problems characterized by high P\'eclet numbers are far more complex and may exhibit instabilities even with classical high fidelity numerical approximations, such as finite element or finite difference method. To deal with this issue we have to resort to some stabilization techniques \cite{brooks, quarteroni_valli}, such as SUPG stabilization. A similar stabilization needs to be accounted for also at the reduced order level, resulting in a stabilized version of the RB algorithm \cite{pacciarini_a, pacciarini_c, pacciarini_b}. In particular, in these works it was shown that a double stabilization in \emph{Offline} and \emph{Online} stage was necessary to obtain an accurate approximation. Nevertheless, stabilizations in \emph{Online} phase can be a bothersome computational cost that may damage the efficiency of the method (for example in \textit{many-query} context), while in some other situation an \emph{Offline--only} stabilized method can be preferred. Stabilization of problems characterized by strong convection effects is an active topic of research in the model order reduction community, see e.g.\ \cite{Akhtar,shafqat,Baiges,Enrique,Traian1,Traian3,Traian2,Lorenzi,maday, pacciarini_a, pacciarini_c,Giovanni} for several different proposed methods with applications in heat transfer and computational fluid dynamics.
\par When dealing with stochastic equations, i.e., with random input parameters, we can modify the RB method, according to probability laws that rule our parameters. In this direction, the wRB (weighted reduced basis) method \cite{peng_art} wants to exploit all the information that random variables give us (a review is provided in \cite{peng_A}). The main novelty of the papers are (i) the synergy of wRB with a stabilized formulation, suitable for stochastic advection dominated problems, and the resulting (ii) capability to enable adaptive toggling of the stabilization depending on the stochastic P\'eclet number. In particular, we will apply the weighted method to stabilized reduced basis strategies and prove the accuracy of the combined method.
%\par Finally, we will study a way to reduce computational costs, combining the \emph{Offline--only} and the \emph{Offline--Online }stabilized reduced basis methods, guaranteeing an error estimation over this procedure. 
%\par
Throughout the work we will test these methods on some steady and time--dependent problems. 

The outline of the manuscript is as follows. In section \ref{stab_RB} we will briefly introduce elliptic coercive parametrized PDEs, their associate RB method, some classical stabilization methods for FE approximation of advection dominated problems; then we will study two reduced basis stabilization methods by testing them on some examples.
We will consider next stochastic partial differential equations; we will present in section \ref{stochastic} the weighted RB method and we will combine it with proper stabilization techniques. Moreover, we will provide a method that selectively enables stabilization to optimize computational costs.
In section \ref{parabolic_RB} we will extend these ideas to parabolic problems, by introducing the general weighted RB method for these problems, combining it with a suitable stabilization technique (based on stabilization for the FE approximation of advection dominated parabolic problems), and testing it on few examples. Finally, section \ref{conclusions} will provide some conclusions and future perspectives.

\section{Stabilized reduced basis method for deterministic elliptic equations}
\label{section stabilized}\label{stab_RB}
\subsection{A brief introduction to reduced basis method}
\label{RB}
The reduced basis (RB) method is a reduced order modelling (ROM) technique which provides rapid and reliable solutions for parametrized partial differential equations (PPDEs) \cite{reduced_basis}, in which the parameters can be either physical or geometrical, deterministic or stochastic.
\par
The need to solve this kind of problems arises in many engineering applications, in which the evaluation of some \emph{output} quantities is required. These \emph{outputs} are often functionals of the solution of a PDE, which can in turn depend on some \emph{input} parameters. The aim of the RB method is to provide a very fast computation of this \emph{input-output} evaluation and so it turns out to be very useful especially in \emph{real-time} or \emph{many-query} contexts.
\par Roughly speaking, given a value of the parameter, the (Lagrange) RB method consists in a Galerkin projection of the continuous solution on a particular subspace of a high-fidelity approximation space, e.g. a finite element (FE) space with a large number of degrees of freedom. This subspace is the one spanned by some pre-computed high-fidelity global solutions (snapshots) of the continuous parametrized problem, corresponding to some properly chosen values of the parameter.
\par For a complete presentation of the reduced basis method we refer to \cite{reduced_basis}, now we just recall its main features in order to
introduce some notations.
\subsubsection{The continuous problem}
Let $\bmu$ belong to the parameter domain $\mathcal{D}\subset \R^p$, $p \in \mathbb{N}$. Let $\Omega$ be a regular bounded open subset of $\R^d$, $d=1,2,3$, and $X$ a suitable Hilbert space.
For any $\bmu\in \mathcal{D}$, let $a(\cdot,\cdot ;\bmu):X \times X\to \R$ be a bilinear form and let $F(\cdot; \bmu):X\to \R$ be a linear functional. 
As we will focus on advection--diffusion equations, that are second order elliptic PDE, the space $X$ will be such that $H^1_0(\Omega)\subset X \subset H^1(\Omega)$. Formally, our problem can be written as follows:
\begin{equation} \label{eq:original_problem}
\begin{split}
&\text{for any }\bmu \in \mathcal{D}, \text{ find }u(\bmu)\in X:\\ 
&a(u(\bmu),v;\bmu)=F(v;\bmu),\quad \forall v \in \mathit{X}.
\end{split}
\end{equation}
We require $a$ to be coercive and continuous, i.e., respectively:
\begin{equation}\label{coercivity}
\exists ~\alpha_0 \text{  s.t.  } \alpha_0 \leq \alpha (\bmu)= \inf_{v \in X} \frac{a(v,v;\bmu )}{||v||_X^2}, \quad \forall \bmu \in \mathcal{D},
\end{equation}
and
\begin{equation}\label{continuity}
+\infty > \gamma (\bmu) = \sup_{v\in X} \sup_{w \in X} \frac{|a(v,w;\bmu)|}{||v||_X||w||_X}, \quad \forall \bmu \in \mathcal{D}.
\end{equation}
For the sake of online efficiency, we assume an \emph{affine} dependence of $a$ on the parameter $\bmu$, i.e.\ we assume that
\begin{equation}\label{affine_a}
a(v,w;\bmu)=\sum_{q=1}^{Q_a} \Theta_a^q(\bmu )a^q(v,w), \quad \forall \bmu \in \mathcal{D}.
\end{equation}
Here, $\Theta_a^q(\bmu ):\mathcal{D}\rightarrow\mathbb{R},$ $q=1,\dots , Q_a$, are smooth functions, while $a^q:X \times X \rightarrow \mathbb{R},$ $q=1,\dots,Q_a$, are $\bmu$-independent continuous bilinear forms. %This assumption will turn out to be crucial for performing the \emph{Offline--Online} decoupling of the computation \cite{reduced_basis}.

In a similar way, we assume that also the functional $F$ is continuous and depends \textquotedblleft affinely\textquotedblright\, on parameters:

\begin{equation}\label{affine_F}
F(v;\bmu)=\sum_{q=1}^{Q_F} \Theta_F^q(\bmu )F^q (v),\quad \forall \bmu \in \mathcal{D},
\end{equation}
where, also in this case, $\Theta_F^q(\bmu ): \mathcal{D} \rightarrow \mathbb{R} ,$ $q=1,\dots , Q_F$, are smooth functions, while $F^q:X \rightarrow \mathbb{R},$ $q=1,\dots,Q_F$, are $\bmu$-independent continuous linear functionals.
\par
Let $X^\mathcal{N} \subset X$ be a conforming finite element space with $\mathcal{N}$ degrees of freedom, we can now set the \emph{truth }approximation of the problem \eqref{eq:original_problem}:
\begin{equation}\label{problem_FE}
\begin{split}
&\text{for any }\bmu \in \mathcal{D}, \text{ find }\uFE\in\mathit{X}^\mathcal{N} ~\text{s.t.}\\
&a(\uFE,v^\mathcal{N};\bmu)=F(v^\mathcal{N};\bmu), \quad \forall v^\mathcal{N} \in \mathit{X}^\mathcal{N}.
\end{split}
\end{equation}
As we are considering the conforming FE case, conditions similar to \eqref{coercivity} and \eqref{continuity} are fulfilled by restriction. More precisely, as regards the coercivity of the restriction of $a$ to $X^\mathcal{N}\times X^\mathcal{N}$, we define:
\begin{equation}\label{coercivity_FE}
\alpha ^\mathcal{N}(\bmu):= \inf_{v^\mathcal{N} \in X^\mathcal{N}} \frac{a(v^\mathcal{N},v^\mathcal{N};\bmu )}{||v^\mathcal{N}||_X^2}, \quad \forall \bmu \in \mathcal{D}
\end{equation}
and, as we are considering a restriction, it easily follows that $
\alpha (\bmu) \leq \alpha ^\mathcal{N} (\bmu), \quad \forall \bmu \in \mathcal{D}$. Similarly, for the continuity, we can define
\begin{equation}\label{continuity_FE}
+\infty > \gamma^\mathcal{N} (\bmu) = \sup_{v^\mathcal{N}\in X^\mathcal{N}} \sup_{w^\mathcal{N} \in X^\mathcal{N}} \frac{|a(v^\mathcal{N},w^\mathcal{N};\bmu)|}{||v^\mathcal{N}||_X||w^\mathcal{N}||_X}, \quad \forall \bmu \in \mathcal{D}.
\end{equation}

%and it holds that $ \gamma (\bmu) \geq \alpha ^\mathcal{N} (\bmu) \quad \forall \bmu \in \mathcal{D}.$ \\
As we have already mentioned, also the domain of the equation can depend on the parameter. In this case we need to map the parametric domain $\Omega_p(\bmu)$ onto a reference one denoted with $\Omega$, via suitable parameter--dependent transformation $T(\cdot;\bmu):\Omega \to \Omega_p(\bmu)$, see \cite{articolo_francesco, reduced_basis, lassila_free_form, manzoni_36}. This allows to track back on the reference domain $\Omega$ all the involved bilinear and linear forms, so that \eqref{affine_a} and \eqref{affine_F} are defined on a common reference domain $\Omega$. In this work we used only affine mappings \cite{reduced_basis, manzoni_36} that allow to easily recover the affinity assumptions \eqref{affine_a} and \eqref{affine_F}. In \cite{manzoni_36,reduced_50} it is possible to find, in particular, a detailed treatment of the advection--diffusion operators.

\subsubsection{The reduced basis method: main features}
\label{sec:rbdet}
Let us suppose that we are given a problem in the form \eqref{eq:original_problem} and its \emph{truth }approximation \eqref{problem_FE}. We recall that the dimension of the finite element space $X^\mathcal{N}$ is $\mathcal{N}$. Given an integer $N\ll\mathcal{N}$, suppose that we are given a set of $N$ suitable parameter values, $S_N=\{ \bmu^1,\dots, \bmu^N\}$: this allows us to define the \emph{reduced basis space} as $X^\mathcal{N}_N=\text{span}\{u^\mathcal{N}(\bmu^n) : 1\leq n\leq N \}$. To be more precise, a Gram-Schmidt orthonormalization process on $\{u^\mathcal{N}(\bmu^n) : 1\leq n\leq N \}$ is usually carried out for the sake of numerical stability, and the resulting orthonormal functions are considered as bases of the reduced space \cite{reduced_basis, reduced_50}.
\par Given a value $\bmu\in\mathcal{D}$, we define the RB solution $\uRB$ such that:
\begin{equation}\label{reduced_equation}
a(u^\mathcal{N}_N (\bmu), v_N; \bmu)=F(v_N;\bmu) \quad \forall v_N\in X_N^\mathcal{N}.
\end{equation}
\par
Recalling that $N \ll \mathcal{N}$, we emphasize the fact that to find the RB solution we need just to solve a $N \times N$ linear system, instead of the $\mathcal{N} \times \mathcal{N}$ one of the FE method.
Moreover, we can also guarantee that the error for a parameter $\bmu \in \mathcal{D}$ is bounded by an error estimator $\Delta_N(\bmu)$:
\begin{equation}\label{def:error_estimator}
|||\uFE -u^\mathcal{N}_N (\bmu) |||_{\bmu} \leq \Delta_N(\bmu) \quad \forall \bmu \in \mathcal{D},
\end{equation}
where $|||\cdot |||_{\bmu}$ is the norm induced by the symmetric part $a_S(\cdot, \cdot;\bmu)$ of the bilinear form $a(\cdot, \cdot;\bmu)$.
The error estimator is defined as $\Delta_N(\bmu):= \frac{||\hat{r}(\bmu)||_X}{\sqrt{\alpha_{LB}(\bmu)}} $, where $\hat{r}$ is the Riesz representor for the functional $r(v^\mathcal{N},\bmu)= F(v^\mathcal{N}; \bmu) - a (u^\mathcal{N}_N(\bmu), v^\mathcal{N}; \bmu) $, $||\cdot ||_X$ is the norm associated to the scalar product in $X$ and $\alpha_{LB}(\bmu)$ is a lower bound for the coercivity constant $\alpha(\bmu)$, possibly dependent on $\bmu\in \mathcal{D}$.
\par The set $S_N$ is built in the \emph{Offline} stage using a Greedy algorithm on a training set $\Xi_{train}$ that spans $\mathcal{D}$ \cite{reduced_basis, reduced_50}. It is an iterative method that, at each step, chooses the parameter value which maximizes the a posteriori error estimator $\bmu\mapsto \Delta_N(\bmu)$ in the training set.
The algorithm stops when a prescribed tolerance $\varepsilon^*_{tol}$ is reached, that is when $\Delta_N(\bmu)\leq  \varepsilon^*_{tol}$ for each parameter value $\bmu$ in the training set $\Xi_{train}\subset \mathcal{D}$. We assume in this section that $\Xi_{train}$ is a collection of randomly selected parameter values according to an uniform distribution.
The error estimator $\Delta_N$ is sharp, in order to avoid an unnecessarily high dimension $N$ for the reduced basis space. Moreover, it must be computationally inexpensive in order to speed up the Greedy algorithm (within which it is computed many times) and to allow the certification of the RB solution during the \emph{Online} stage.
\par We want to point out that all the expensive computations (i.e. those whose costs depend on the FE space dimension $\mathcal{N}$) are performed during the \emph{Offline} stage.
Indeed, the affinity assumptions \eqref{affine_a} and \eqref{affine_F} are crucial for the \emph{Offline---Online} decoupling, as it is extensively shown in \cite{reduced_basis,reduced_50}. The affinity assumptions allow the storage, during the \emph{Offline} stage, of the matrices corresponding to the parameter independent forms $a_q, q=1,\dots,Q_a$, restricted to $X^\mathcal{N}_N$. Thanks to this fact, during the \emph{Online} stage the assembly of the reduced basis system only consists in a linear combination of these precomputed matrices. A similar strategy can  also be applied to the computation of the error estimator \cite{reduced_basis, reduced_50}.
Indeed, thanks to the affine decomposition of $F$ \eqref{affine_F} and $a$ \eqref{affine_a}, $\hat{r}$ can be computed in an \emph{Online} phase, with a complexity that only depends on $N$ but not on $\mathcal{N}$ \cite{reduced_basis}. Also the $\alpha_{LB}(\bmu)$ can be efficiently computed in an \emph{Online} phase, thanks to suitable algorithms such as the successive constraint method \cite{reduced_basis, huynh_scm}. Therefore, at each step of the Greedy algorithm, the error estimator $\Delta_N(\bmu)$ can be efficiently evaluated (with computational complexity independent from $\mathcal{N}$) for any element in the training set, rather than relying on the computation of the error $|||\uFE -u^\mathcal{N}_N (\bmu) |||_{\bmu}$ (which would require an expensive \emph{truth} solve for all parameters in the training set, such as in a proper orthogonal decomposition basis generation). 
If affinity assumptions are not fulfilled, it turns out to be necessary to use an interpolation strategy (e.g. empirical interpolation method (EIM) \cite{barrault,eftang_EIM}) in order to recover them. A weighted version of EIM is provided in \cite{peng_AAA}.
\subsection{Stabilized reduced basis methods}
\par
The main goal of this section is to design an efficient stabilization procedure for the RB method. More specifically, we will make a comparison between an \emph{Offline--Online} stabilized method and an \emph{Offline--only} stabilized one as done in \cite{pacciarini_a}. We want to approximate the solution of a parametric advection--diffusion problem:

\begin{equation}\label{adv-diff_equation}
-\varepsilon \Delta u + \bbe \cdot \nabla u = f\quad \text{in } \Omega \subset \mathbb{R}^d
\end{equation}
given a parameter value $\bmu\in\mathcal{D}$ and suitable Dirichlet, Neumann or mixed boundary conditions. Here $\varepsilon = \varepsilon(\bmu): \Omega \to [0,+\infty)$ is a parametrized diffusion coefficient, while $\bbe = \bbe(\bmu): \Omega\to \R^d$ is a parametrized advection field such that div$(\bbe)=0$.

%It is well known from the general theory of the numerical approximation of advection--diffusion equations, that the FE solution of such equations can show significant instability phenomena when the advective term dominates the diffusive one \cite{quarteroni_valli}. 
Let $\mathcal{T}_h$ be a triangulation of $\Omega$ and let $K$ be an element of $\mathcal{T}_h$. We say that a problem is advection dominated in $K$ if the following condition holds:
\begin{equation}\label{peclet_adv_dom}
\pec := \frac{|\bbe (x) | h_k}{2 \varepsilon(x)}>1\quad \forall x \in K,
\end{equation}
where $h_K$ is the diameter of K.
It is very well known from literature (e.g. \cite{quarteroni_valli}) that the FE approximation of advection dominated problems can show significant instability phenomena, e.g. spurious oscillations near the boundary layers. 
Several recipes have been proposed to fix these issues. We choose to resort to a strongly consistent stabilization method: the Streamline/Upwind Petrov--Galerkin (SUPG) \cite{brooks, hughes_brooks, johnson1, johnson2}. The main idea of stabilization techniques is to add artificial diffusion to equation \eqref{adv-diff_equation}. To increase the accuracy of the resulting solution, SUPG adds diffusion only in the streamline direction, and not everywhere as in a purely artificial diffusion scheme. Moreover, the resulting method is strongly consistent with the continuous PDE and, provided that the stabilization coefficients are properly chosen, retains the same order of accuracy as the underlying discretization scheme.
For a detailed presentation of the stabilization method for the FE approximation of advection dominated problems, we refer to \cite{hughes_brooks, quarteroni_valli}.
\par 
Let us now explain the basic ideas of the two RB stabilization methods mentioned before. As regards the \emph{Offline--Online} stabilized method, the choice of the name reveals that the Galerkin projections are performed, in both \emph{Offline} and \emph{Online} stage, with respect to the SUPG stabilized bilinear form \cite{brooks, quarteroni_valli}, that is 
\begin{align}
\label{eq:lhs} a_{stab}(w^\mathcal{N},v^\mathcal{N};\bmu) &= a(w^\mathcal{N},v^\mathcal{N};\bmu) + s(w^\mathcal{N},v^\mathcal{N};\bmu)\\
F_{stab}(v^\mathcal{N};\bmu) & = F(v^\mathcal{N};\bmu) + r(v^\mathcal{N};\bmu)\notag\\
\label{no_stabilization_lhs} a(w^\mathcal{N},v^\mathcal{N};\bmu) &= \int_{\Omega} \varepsilon(\bmu) \nabla w^\mathcal{N} \cdot \nabla v^\mathcal{N}+(\bbe(\bmu) \cdot \nabla w^\mathcal{N})v^\mathcal{N}\\
F(v^\mathcal{N};\bmu) & = \int_{\Omega} f v^\mathcal{N}\notag\\
\label{stabilization_lhs} s(w^\mathcal{N},v^\mathcal{N};\bmu) &= \sum _{K\in \mathcal{T}_h} \delta_K  \int_K L w^\mathcal{N} \frac{h_K}{|\bbe(\bmu) |}L_{SS}v^\mathcal{N}\\
r(v^\mathcal{N};\bmu) &= \sum _{K\in \mathcal{T}_h} \delta_K \int_K f \frac{h_K}{|\bbe(\bmu) |}L_{SS}v^\mathcal{N}
\end{align}
where $w^\mathcal{N},v^\mathcal{N}$ chosen in a suitable piecewise polynomial space $X^\mathcal{N}$. In \eqref{stabilization_lhs} $L$ is the parameter dependent advection--diffusion operator, that is $Lv^\mathcal{N}=\varepsilon \Delta v^\mathcal{N} + \bbe \cdot \nabla v^\mathcal{N}$, which can be split into its symmetric part $L_Su^\mathcal{N}=-\varepsilon \Delta u^\mathcal{N}$ and its skew--symmetric part $L_{SS}u^\mathcal{N}= \bbe \cdot \nabla u^\mathcal{N}$. Moreover, $h_K$ denotes the diameter of the element $K$, while $\delta_K$ is a positive real number which may depend on $K$ through the parameter $\bmu$ (but not directly on $h_K$).
\par
%The bilinear form $a_{stab}$ is coercive, so we can apply the already developed theory in order to use the reduced basis method. 
In contrast, in the \emph{Offline--only} stabilized method we use the stabilized form \eqref{eq:lhs} only during the \emph{Offline} stage, while during the \emph{Online} stage we project with respect to the standard advection--diffusion bilinear form \eqref{no_stabilization_lhs}. 
%To motivate this choice, we recall that the RB approximation is actually a linear combination of particularly chosen \emph{truth} solutions of the problem. Thus, if we consider stable \emph{truth} approximations when building the reduced basis, it seems reasonable to expect RB approximations which does not show instability. Moreover, note that the \emph{Offline--only }method is strongly consistent with respect to the continuous advection--diffusion problem.
%\par
An advantage of using the \emph{Offline--only }stabilized method would be a certain reduction of the online computational effort in the assembly of the reduced linear system, that could be also significant if the number of affine stabilization terms is very high. Among possible disadvantages, we mention the inconsistency between the offline and online bilinear forms.
\par
We will start from the study of some test problems, which we will keep as prototypes for each further extension that will be carried out in the next sections. The first one is a \PGP problem \cite{incropera,manzoni_45,pacciarini_a}, while the second is a parametrized internal layer problem \cite{pacciarini_a}. From here on, we will explicitly write the FE space dimension $\mathcal{N}$ only when it will be strictly necessary.

\subsubsection{Numerical test: Poiseuille--Graetz problem (PG)}\label{Graetz}

\begin{figure}[h]
\centering
\includegraphics[trim=7cm 14cm 7cm 15cm, clip=true, scale=0.28]{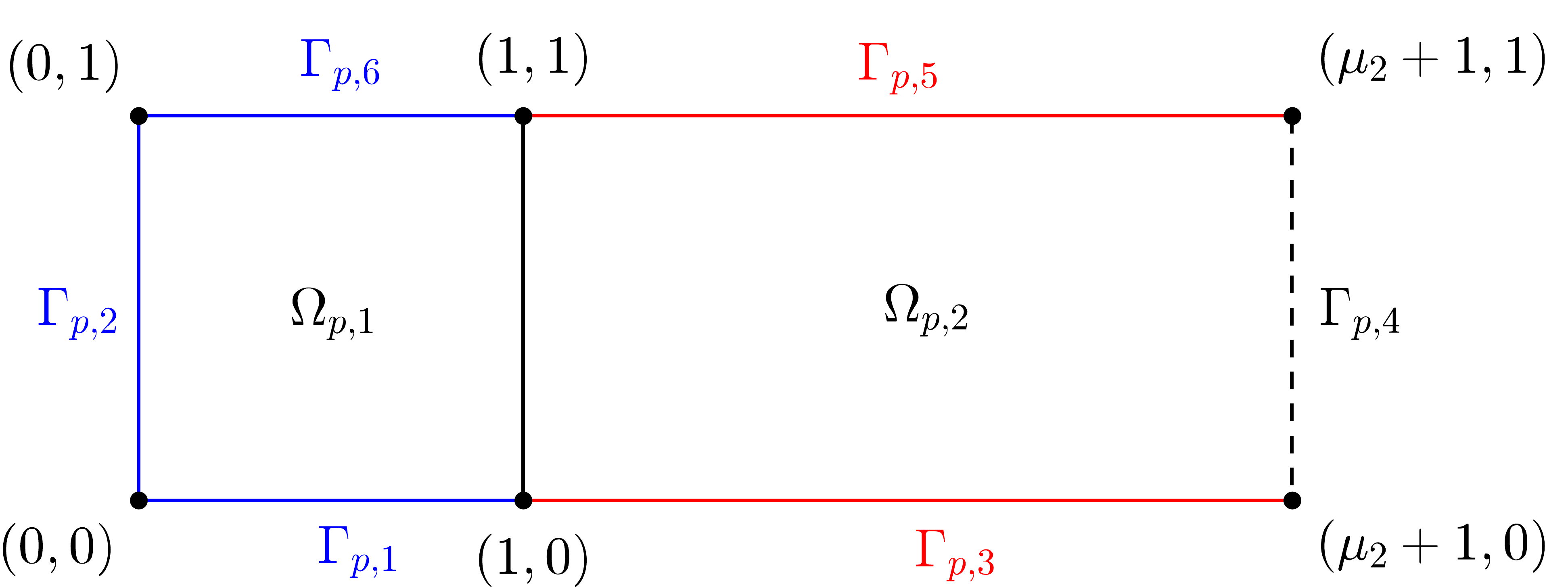}
\caption{Geometry of \PGP problem. Parametrized domain. Boundary conditions: homogeneous Dirichlet on blue sides, $u=1$ on red sides, homogeneous Neumann on the dashed side}\label{graetz_geometry}
\end{figure}

We consider a \PGP problem where we have two parameters: one physical (the inverse of diffusivity coefficient $\mu_1$, which is proportional to the P\'eclet number) and one geometrical (the length of the domain being equal to $1+\mu_2$). The \PGP problem deals with steady forced heat convection (advective phenomenon) combined with heat conduction (diffusive phenomenon) in a duct with walls at different temperature. Let us define $\bmu =(\mu_1, \mu_2)$ with both $\mu_1$ and $\mu_2$ positive, real numbers. Let $\Omega_p(\bmu)$ be the rectangle $(0,1+\mu_2)\times (0,1)$ in $\mathbb{R}^2$. The domain is shown in figure \ref{graetz_geometry}.
The problem is to find a solution $u(\bmu)$, representing the temperature distribution, such that:
\begin{equation}\label{graetz_strong}
\begin{cases}
-\frac{1}{\mu_1} \Delta u(\bmu) +4 y(1-y) \partial_x u(\bmu) = 0 & \text{in } \Omega_p(\bmu)\\
u(\bmu)=0 & \text{on } \Gamma_{p,1}(\bmu)\cup  \Gamma_{p,2}(\bmu)\cup  \Gamma_{p,6}(\bmu)\\ 
u(\bmu)=1 &\text{on } \Gamma_{p,3}(\bmu)\cup  \Gamma_{p,5}(\bmu)  \\
\frac{\partial u}{\partial \nu }=0 & \text{on } \Gamma_{p,4}(\bmu).
\end{cases}
\end{equation}
\par
We set the reference domain as $\Omega = (0,2)\times (0,1)$, and subdivide it in $\Omega^1= (0,1)\times (0,1)$ and $\Omega^2= (1,2)\times (0,1)$. The affine transformation that maps the reference domain into the parametrized one is:
\begin{align}\label{graetz_trasform_1}
&T^1(\bmu):\Omega^1 \to \Omega_{p,1}(\bmu) \subset \mathbb{R}^2 \hfill &T^2(\bmu):\Omega^2 \to \Omega_{p,2}(\bmu) \subset \mathbb{R}^2\\
&T^1\left(\begin{pmatrix}
x\\ y \end{pmatrix};\bmu\right)=\begin{pmatrix} x\\y \end{pmatrix} \hfill &T^2\left(\begin{pmatrix}
x\\y
\end{pmatrix};\bmu\right)=\begin{pmatrix}
\mu_2 x\\y
\end{pmatrix}  +\begin{pmatrix}
1-\mu_2\\0
\end{pmatrix} .
\end{align}
and define the continuous one--to--one transformation $T(\bmu)$ by gluing together these two transformations.

Let us now define a mesh $\mathcal{T}_h$ on the reference domain $\Omega$ and let us call $\mathcal{T}_h^1$ and $\mathcal{T}_h^2$ the restrictions $\mathcal{T}_h$ to $\Omega_1$ and $\Omega_2$, respectively. We use $\mathbb{P}^1$ FE discretization during the offline stage.
%We can also define a mesh on $\Omega_p(\bmu)$ just by taking the image of $\mathcal{T}_h$ through the transformation $T(\cdot,\bmu)$. 
Hence, the corresponding bilinear forms $a(\cdot, \cdot; \bmu)$ and $s(\cdot, \cdot; \bmu)$ are
\begin{equation}\label{Graetz_supg_lhs_a}
\begin{split}
a(u^\mathcal{N},v^\mathcal{N};\bmu):=&\int_{\Omega^1} \frac{1}{\mu_1}\nabla u^\mathcal{N} \nabla v^\mathcal{N} + 4y(1-y) \partial_x u^\mathcal{N} v^\mathcal{N}+\\
&+\int_{\Omega^2} \frac{1}{\mu_1\mu_2} \partial_x u^\mathcal{N} \partial_xv^\mathcal{N} + \frac{\mu_2}{\mu_1} \partial_x u^\mathcal{N} \partial_y v^\mathcal{N} + 4\mu_2 y(1-y) \partial_x u^\mathcal{N} v^\mathcal{N}
\end{split}
\end{equation}
and
\begin{equation}\label{Graetz_supg_lhs_s}
\begin{split}
s(u^\mathcal{N},v^\mathcal{N};\bmu):=\sum_{K\in\mathcal{T}_h^1} h_K \int_{K}( 4y(1-y) \partial_x u^\mathcal{N}) \partial_x v^\mathcal{N}+\sum_{K\in\mathcal{T}_h^1} \frac{h_K}{\sqrt{\mu_2}} \int_{K} ( 4y(1-y) \partial_x u^\mathcal{N}) \partial_x v^\mathcal{N}.
\end{split}
\end{equation}
The choice of the stabilization coefficient $\delta_{K_p}=\delta_{K_p}(\bmu)=\frac{1}{\sqrt{\mu_2}}$ for $K_p\in \mathcal{T}^2_{h}$ is motivated by the transformation to the reference domain.

\begin{figure}[h]
\centering
\subfigure[$\mu_1 \in (10^4, 10^5 )$ ]{\includegraphics[trim=0cm 0cm 2cm 0cm,clip=true, scale=0.40]{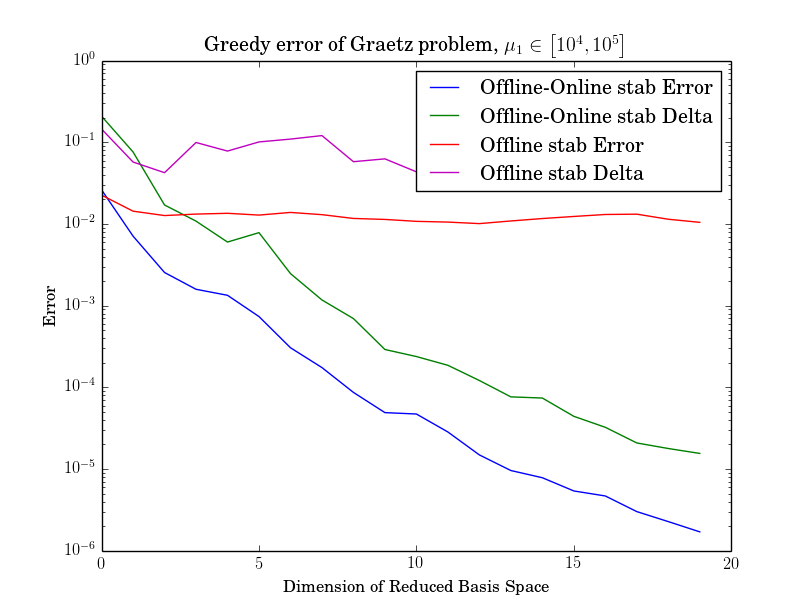} \label{Error_off_offonl_mu_4_5} }
\subfigure[$\mu_1 \in ( 1, 10^6 ) $]{\includegraphics[trim=1.5cm 0cm 2cm 0cm,clip=true, scale=0.40]{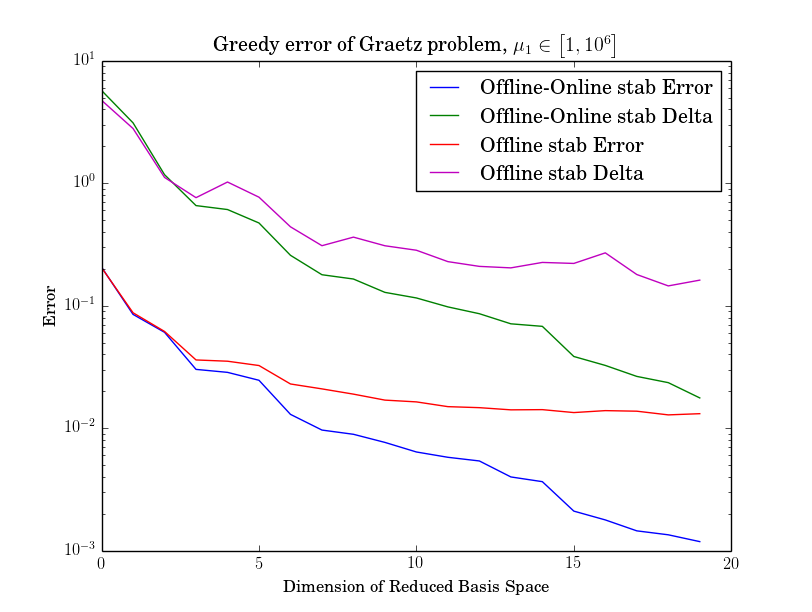} \label{Error_off_offonl_mu_1_6}}
\caption{Error comparison between Offline and Online-Offline stabilization}\label{graetz_error_comp_onl_off}
\end{figure}

We test the performance of the RB approximation for two choices of the parameter space, namely $\mathcal{D}^1= \left[10^4, 10^5 \right] \times \left[ 0.5, 4 \right]$ and $\mathcal{D}^2= \left[1, 10^6\right] \times \left[ 0.5, 4 \right]$. The parameter space $\mathcal{D}^1$ features very large values of $\mu_1$, so that the solution manifold is characterized by solution with steep boundary layers. In contrast, the parameter space $\mathcal{D}^2$ features both small and large values of $\mu_1$, resulting in a richer set of solutions. The range of variation for the geometrical parameter $\mu_2$ is the same in both parameter spaces.
\par
The comparison of \emph{Offline--only} and \emph{Offline--Online} stabilized algorithms is shown in figure \ref{graetz_error_comp_onl_off}, for $\mathcal{D}^1$ (left) and $\mathcal{D}^2$ (right). In each figure, the evolution of the Greedy parameter selection is presented, plotting both the error bound $\max_{\bmu \in \Xi_{train}} \Delta_N(\bmu)$ employed by the RB algorithm and, for comparison, the energy norm error $\max_{\bmu \in \Xi_{train}} |||u^\mathcal{N} (\bmu) - u^\mathcal{N}_N (\bmu)|||_{\bmu}$. For both $\mathcal{D}^1$ and $\mathcal{D}^2$, the Greedy algorithm in the \emph{Online--Offline} case is clearly converging as the RB space enriches its dimension. In contrast, the Greedy algorithm does not converge in the \emph{Offline--only} case, being over $10^{-2}$ for both $\mathcal{D}^1$ and $\mathcal{D}^2$.

\begin{figure}[h]
\centering
\subfigure[Offline-Online stabilized, $\bmu=(10^{4.8},3.3)$]{\includegraphics[trim=0cm 14cm 2cm 9cm,clip=true, scale=0.2]{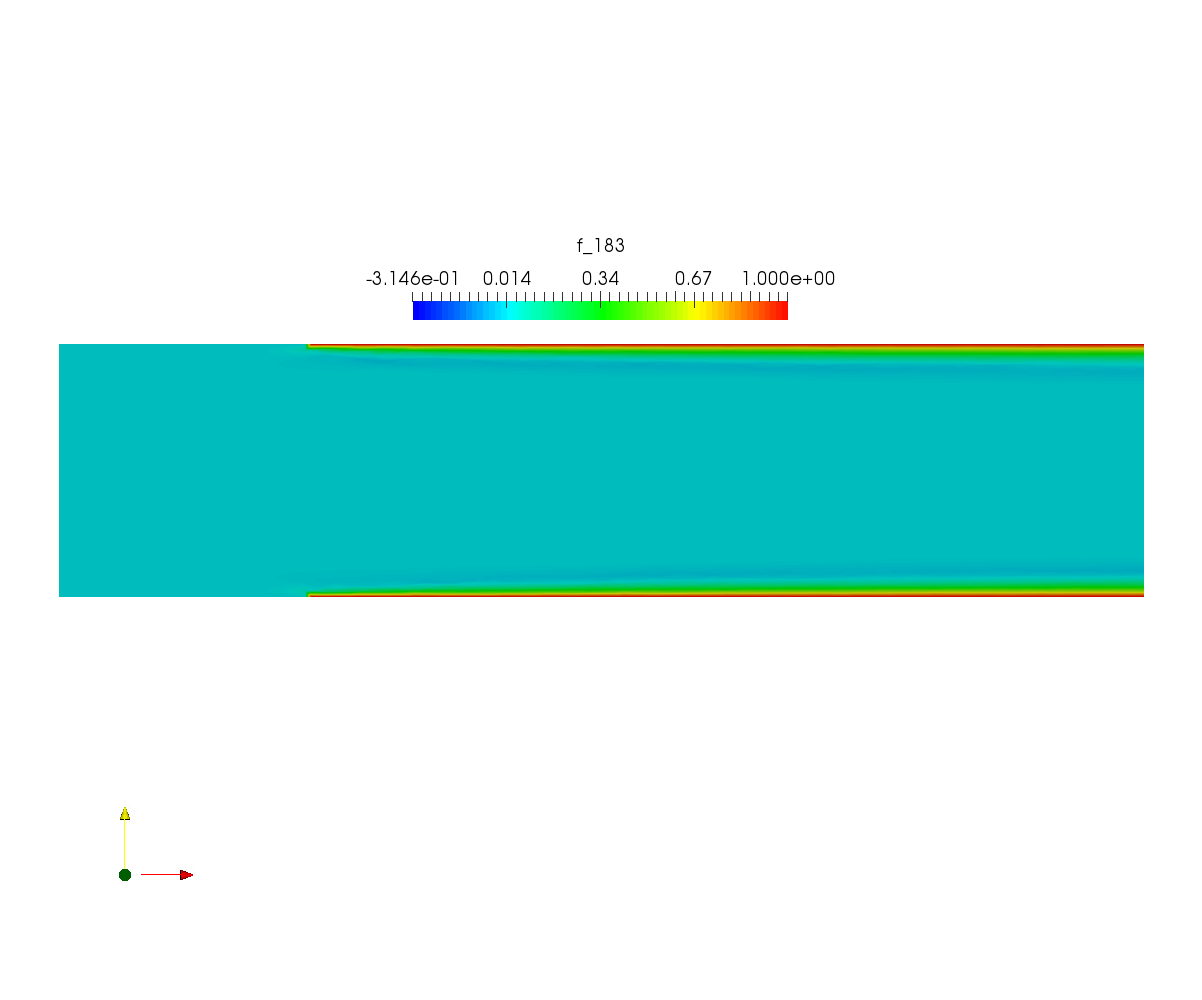} \label{Graetz_stab_RB_off_onl_mu48}}
\subfigure[Zoom on the boundary layer]{\includegraphics[trim=10cm 22cm 30cm 12cm,clip=true, scale=1.6]{Immagini/Graetz/stabilized/RB_mu_33_48_onloff.png} \label{Graetz_stab_RB_off_onl_mu48_zoom}}\\
\subfigure[Offline Stabilized, $\bmu=(10^{4.8},3.3)$]{\includegraphics[trim=0cm 14cm 2cm 9cm,clip=true, scale=0.2]{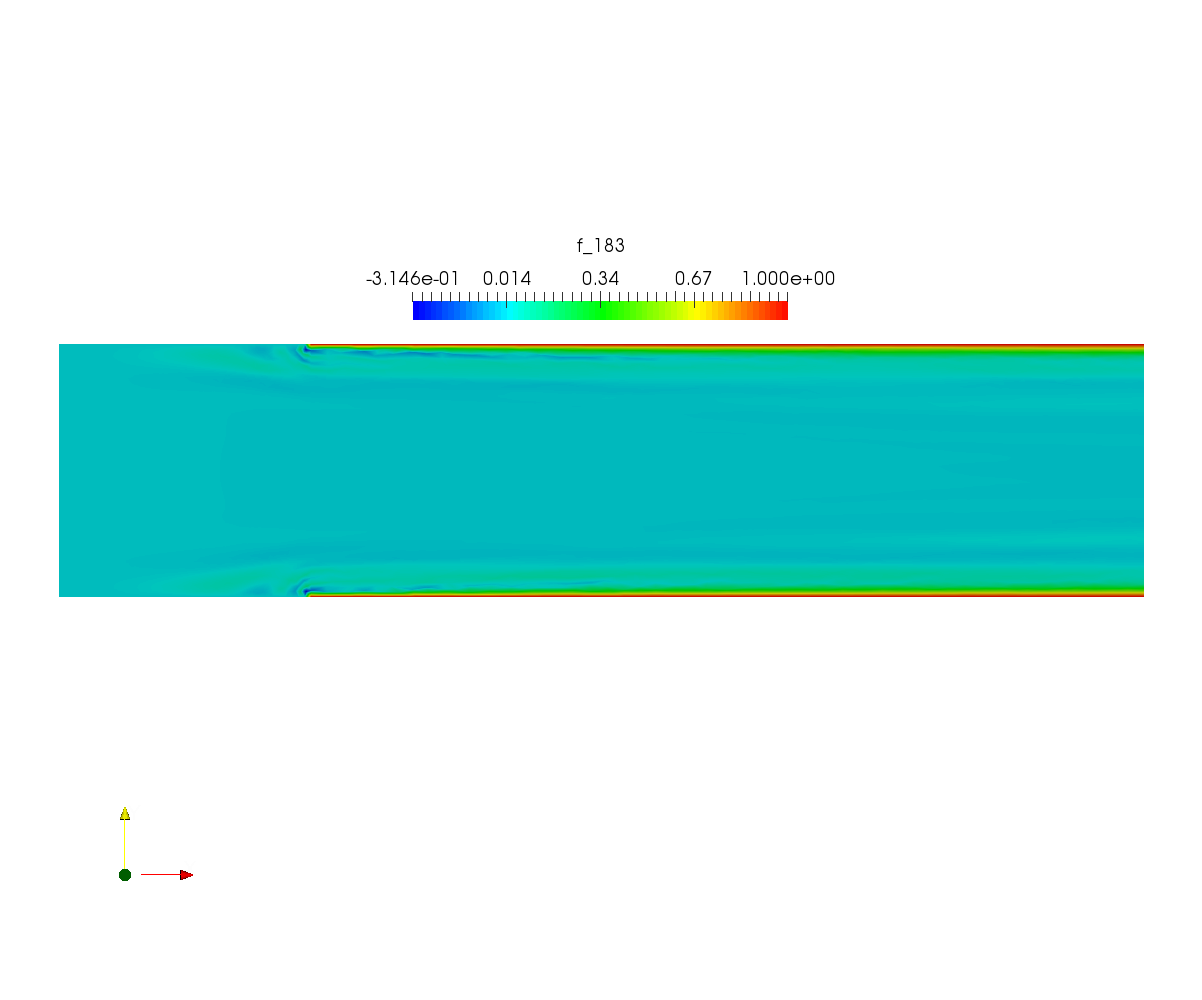} \label{Graetz_stab_RB_off_mu48}}
\subfigure[Zoom on the boundary layer]{\includegraphics[trim=10cm 22cm 30cm 12cm,clip=true, scale=1.6]{Immagini/Graetz/stabilized/RB_mu_33_48_off.png} \label{Graetz_stab_RB_off_mu48_zoom}}\\
\caption{RB solution, stabilized Offline-Online and Offline, $\bmu=(10^{4.8},3.3)$}\label{graetz_stab_RB_10_48}
\end{figure}

We show a representative online solution for both stabilization cases, characterized by large value of P\'eclet number, in figure \ref{graetz_stab_RB_10_48}, obtained for $N=20$. As we can see, the \emph{Offline--Online }stabilized RB solution is showing marked boundary layers, while the \emph{Offline--only} stabilized RB solution still has some noise near the boundary layer and some peaks near discontinuities of solution at top and bottom walls.\\
\par Moreover, if we compare the time used to perform one \textit{truth} solution ($\mathcal{N}=4369$) and a RB one ($N=20$), we can see that the former lasts \textrm{0.0411} seconds, while the stabilized Online RB solution lasts \textrm{0.000512} seconds, on average on a test set. The non-stabilized in the online phase lasts even less time, namely \textrm{0.000151} seconds, even though it is less accurate (see figure \ref{graetz_stab_RB_10_48}).
The further speedup of the non-stabilized version is due to the lower number of affine terms to be assembled online.
Even bigger gains can be observed in the parabolic case in section \ref{parabolic_RB}, or for problems characterized by a large number of affine terms $Q_a$ and $Q_F$.

\begin{comment}
\par
Nevertheless, there exists an error bound for this \emph{Offline--only} error which is sharper than RB one \cite{tesi_paolo, pacciarini_a}. Of course it will be of order of $h_K$ as the stabilization that we have introduced has the same order, and it depends on the tolerance $\varepsilon^*$ of the Greedy algorithm to compute the reduced space $N$. One can prove that the error between stabilized FE solution and \emph{Offline--only }stabilized RB solution is such that \cite{pacciarini_a}:
\begin{equation}\label{error_paccia}
\begin{split}
|||u_N(\bmu)-u^{stab,\mathcal{N}}(\bmu)|||_{\bmu} \leq & h_{max}(\bmu) C(\bmu) ||\bbe \cdot \nabla u^{stab,\mathcal{N}}(\bmu)||_{L^2(\Omega_p(\bmu))}+\\
& +(1+h_{max}(\bmu) C(\bmu)^2 ||\bbe ||_{L^\infty (\Omega_p (\bmu))})\varepsilon^*.
\end{split}
\end{equation}
Here $C(\bmu)$ is a constant used in a norm equivalence \cite{pacciarini_a} and $h_{max}$ is the maximum mesh size.\\
We point out that this bound depends on the $L^2$ norm of the streamline derivative. This means that the \emph{Offline--only} method has better performances when applied to problems in which the strongest variations occur along a direction orthogonal to the advection field. This could happen in the cases in which boundary layers are parallel to the advection field, e.g. the \PGP problem.
\par
We will resort to the \emph{Offline--only }stabilization in the stochastic approach.
\end{comment}
\subsubsection{Numerical test: propagating front in a square (PFS)}\label{section:square}
\begin{figure}[h]
\centering
\includegraphics[trim=21cm 52cm 17cm 15cm,clip=true, scale=0.17]{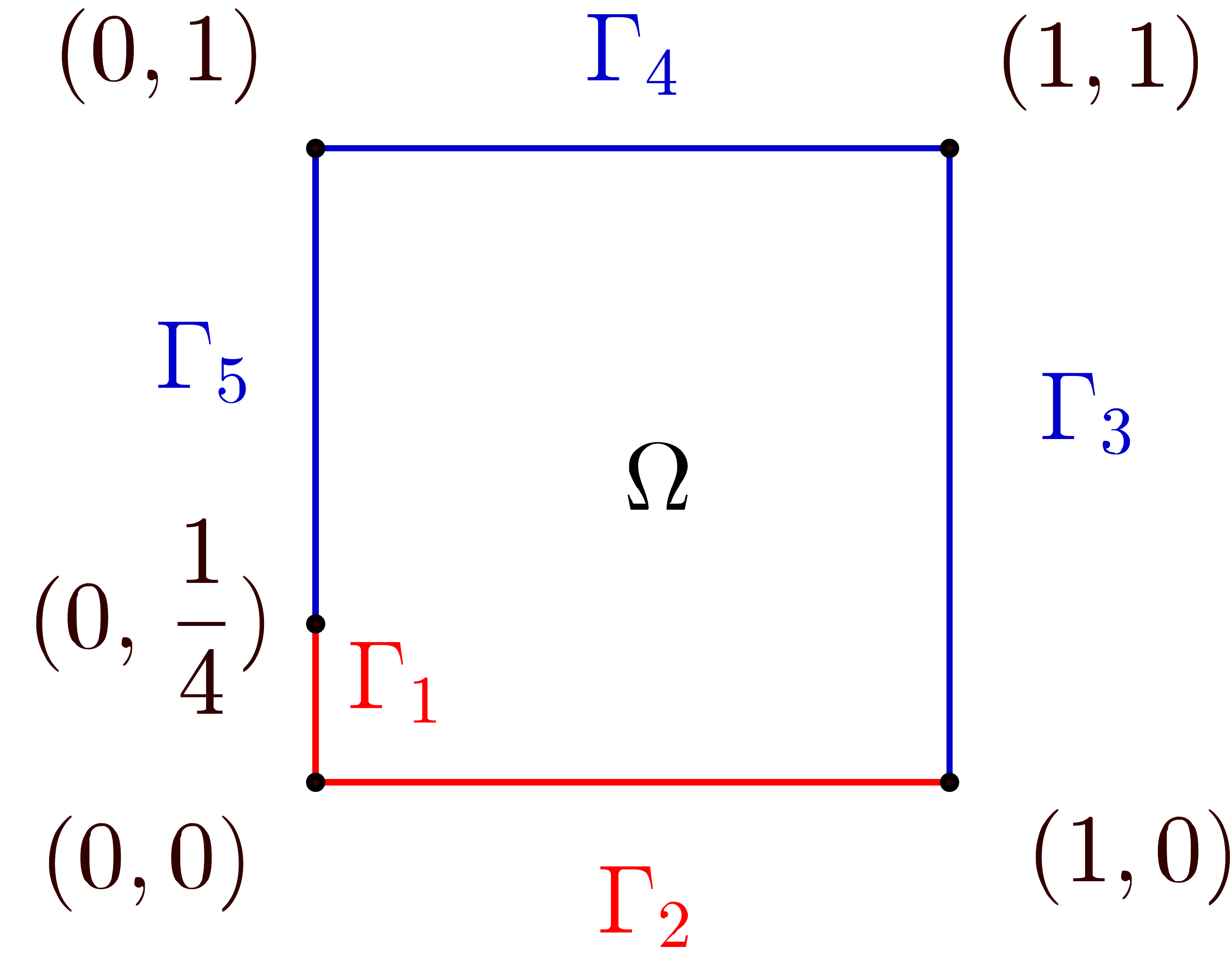}
\caption{Geometry \PFS problem}\label{square_geo}
\end{figure}

In this section we will test the reduced order stabilization method for a second test case where the parameter controls the angle of an internal layer.
The problem we want to study is set over a unit square $\Omega\subset \mathbb{R}^2$, as sketched in figure \ref{square_geo}, it has two parameter $\mu_1,\mu_2\in \mathbb{R}$, and is as follows:
\begin{equation}\label{square_equation}
\begin{cases}
-\frac{1}{\mu_1} \Delta u(\bmu) + (\cos\mu_2, \sin \mu_2 )\cdot \nabla u(\bmu) = 0 & \text{in }\Omega\\
u(\bmu)=1 & \text{on } \Gamma_1 \cup \Gamma_2\\
u(\bmu)=0 & \text{on } \Gamma_3 \cup \Gamma_4 \cup \Gamma_5.
\end{cases}
\end{equation}

\begin{figure}
\centering
\subfigure[$\mu_2=0, \delta_K=2.1$ ]{\includegraphics[trim=0cm 2cm 5cm 0cm,clip=true, scale=0.25]{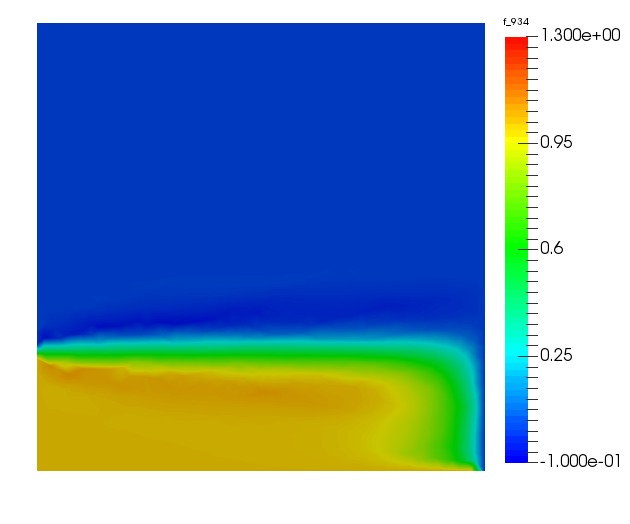} \label{immagine_troopo_stabilizzata} }
\subfigure[$\mu_2=0.8, \delta_K=1.4$ ]{\includegraphics[trim=0cm 2cm 5cm 0cm,clip=true, scale=0.25]{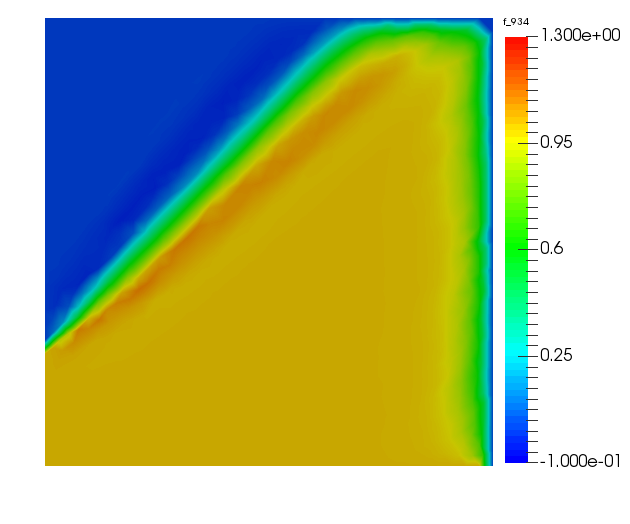} }
\subfigure[$\mu_2=1.2, \delta_K=0.7$ ]{\includegraphics[trim=0cm 2cm 0cm 0cm,clip=true, scale=0.25]{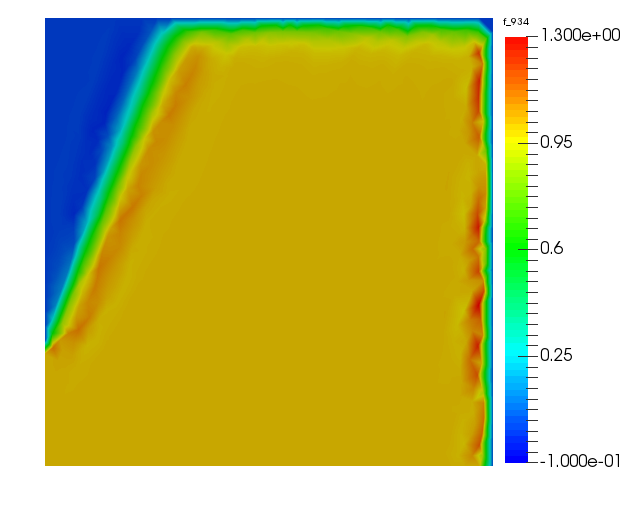} \label{immagine_troppo_poco_stabilizzata} }
\caption{FE solution comparison varying $\delta_K$ and $\mu_2$}\label{square_delta_comparison}
\end{figure}

Let us note that $\mu_1$ is proportional to the P\'eclet number of the advection--diffusion problem, while $\mu_2$ is the angle between the $x$ axis and the direction of the constant advection field. The bilinear form associated to the problem is:
\begin{equation}\label{square_bilinear_lhs}
a(u,v;\bmu)=\int _{\Omega} \frac{1}{\mu_1} \nabla u \cdot \nabla v + (\cos \mu_2 \, \partial_x u + \sin \mu_2 \, \partial _y u)v.
\end{equation}
We introduce again a triangulation $\mathcal{T}_h$ on the domain $\Omega$ and we consider a $\mathbb{P}^1$ discretization. The corresponding stabilization term is
\begin{equation}\label{square_stab_term}
s(u^\mathcal{N},v^\mathcal{N};\bmu) = \sum_{K\in \mathcal{T}^\mathcal{N}} \delta_K \int_{K} (\cos \mu_2, \sin \mu_2)\cdot \nabla u^\mathcal{N} \,\, (\cos\mu_2, \sin \mu_2 ) \cdot \nabla v^\mathcal{N} 
\end{equation}
where $\delta_K$ is manually tuned according to $\mu_2$. 
A few representative FE solutions are shown in figure \ref{square_delta_comparison}. The figure clearly shows that the direction of the advection fields largely affects the solution, which exhibits strong variations in energy norm \cite{tesi_paolo}. For this reason, we test the RB method for two different choices of the parameter space, namely $\mathcal{D}^1 = \left[10^4, 10^5 \right] \times \left[0.5, 1\right]$ and $\mathcal{D}^2 = \left[10^4, 10^5 \right] \times \left[0, 1.57\right]$. Both choices are characterized by dominant advection; moreover, a wider range of angles is considered in $\mathcal{D}^2$ than in $\mathcal{D}^1$, resulting in a richer manifold of solutions.

\begin{figure}[!ht]
\centering
\subfigure[$\mu_2\in (0.5,1)$ ]{\includegraphics[trim=0cm 0cm 1.8cm 0cm,clip=true, scale=0.385]{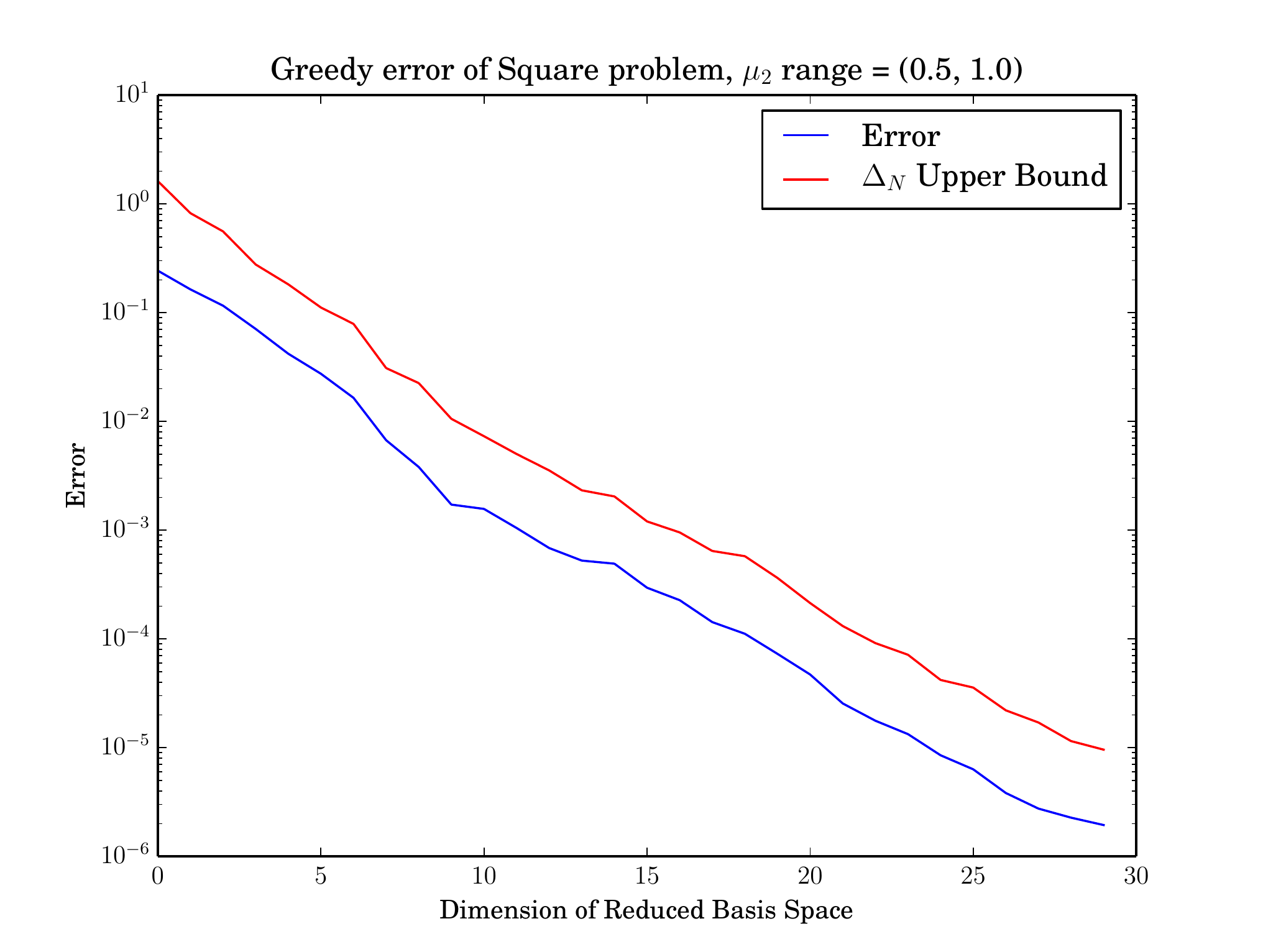} }
\subfigure[$\mu_2\in (0,1.57 )$ ]{\includegraphics[trim=0cm 0cm 1.8cm 0cm,clip=true, scale=0.385]{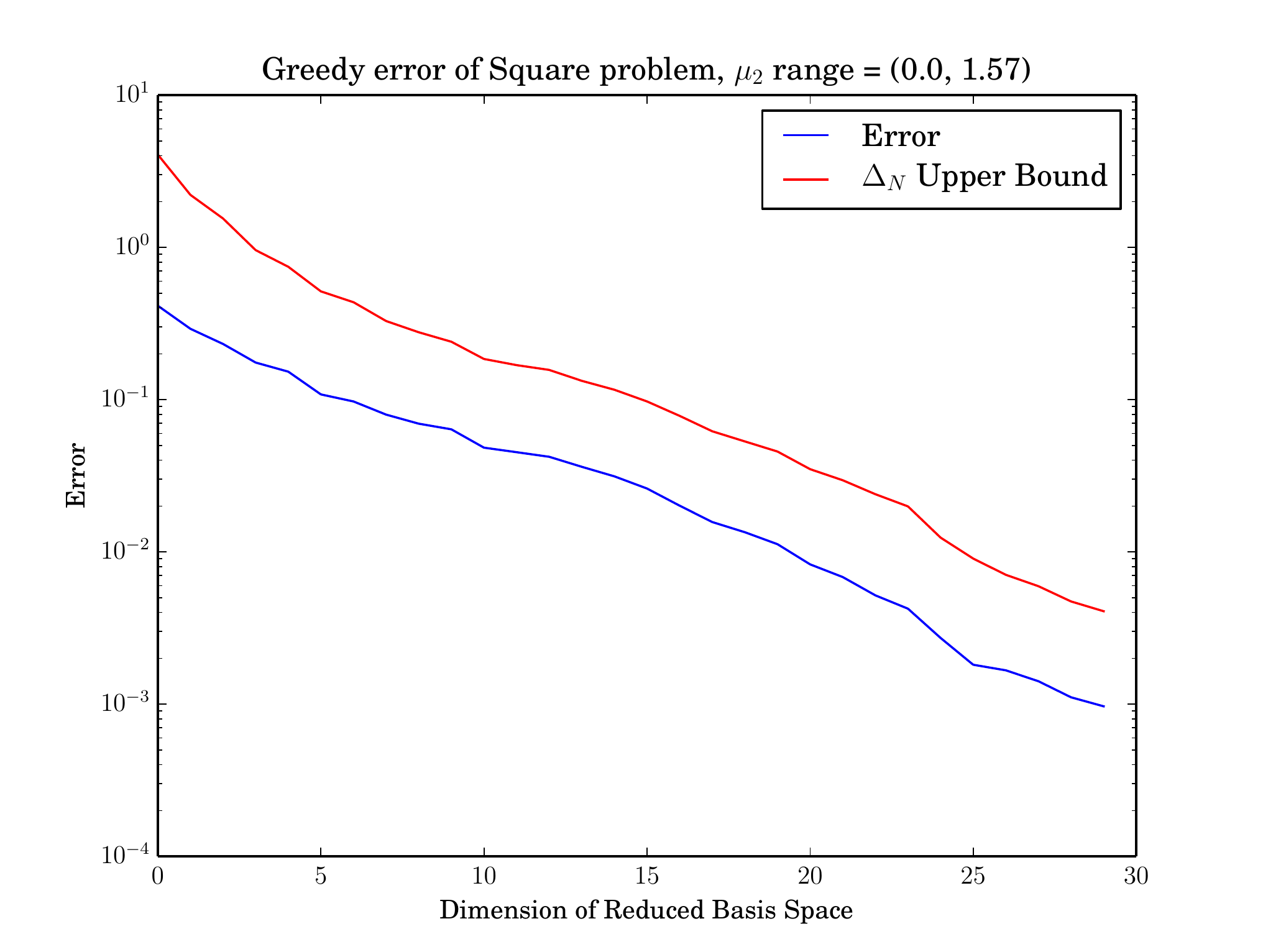} }
\caption{RB error and $\Delta_N$ error bound varying $\mu_2$ range}\label{square_RB_range_comparison}
\end{figure}

The performance of the RB algorithm is shown in figure \ref{square_RB_range_comparison} for $\mathcal{D}^1$ (left) and $\mathcal{D}^2$ (left). Only the \emph{Offline--Online} stabilization case is reported, since the \emph{Offline-only} case gave poor results as in the previous test case.
In both cases the stabilized reduced order method converges, reaching an error around $10^{-6}$ for $\mathcal{D}^1$ and around $10^{-3}$ for $\mathcal{D}^2$. Computational times are: \textrm{0.461346} seconds on average for a \textit{truth} solution ($\mathcal{N}=15626$), \textrm{0.034271} seconds for a RB solution ($N=20$) with online stabilization, and \textrm{0.001862} seconds for a RB solution ($N=20$) without online stabilization.

\section{Stabilized weighted reduced basis algorithm for problems with uncertain parameters}\label{stochastic}
The reduced basis method formulated in section \ref{section stabilized} assumed deterministic parameters; in contrast, for random parameters, a \emph{weighted} reduced basis has been proposed \cite{peng,peng_art} as an extension of the standard reduced basis approach. The main idea of this method is to suitably assign a larger weight to those samples that are more ``important''. 
%We expect a better convergence of this algorithm over a distributed test set.
In this section, we will deal with problems with random distributed parameters and we will compare the weighted method to the standard reduced basis method for advection--diffusion problems with high P\'eclet number. 
%In order to do this, we will use all the stabilization techniques used in section \ref{section stabilized}. 
Moreover, we will also provide an offline/online stabilization approach that can be useful in case when stabilization involves large computations.
%\par Finally, we will provide some numerical examples on problems presented in previous chapters.
%\par Similar methods are studied even for nonlocal diffusion problems. In \cite{gunzburger} we can see both reduced basis methods and uncertainty quantification in sense of random input parameters applied to nonlocal diffusion problems. This tells us that our weighted reduced basis method can be extended to other problems different from PDEs we dealt in this work.
\subsection{A brief introduction to weighted reduced basis method}\label{sec:weighted_RB}
To discuss the weighted reduced basis method \cite{peng_art}, we introduce stochastic partial differential equations. Let $\Omega$ be an open set of $\mathbb{R}^d$ with Lipschitz boundary $\partial \Omega$ and let $H^1_0(\Omega)\subset X \subset H^1(\Omega)$ a functional space. Let $(A,\mathcal{F},P)$ denote a complete probability space, where $A$ is a set of outcomes $\omega\in A$, $\mathcal{F}$ is a $\sigma$-algebra of events and $P:\mathcal{F}\to \left[0,1 \right]$ with $P(A)=1$ is a probability measure \cite{prob}. A real-valued \emph{random variable} is defined as a measurable function $Y:(A,\mathcal{F})\to (\mathbb{R},\mathcal{B})$, being $\mathcal{B}$ the Borel $\sigma$-algebra on $\mathbb{R}$. Let $dF_Y (y)$ denote the distribution measure, i.e., for all $B\subset \mathcal{D}$, $P(F\in B)=\int_B dF_Y(y)$. Provided that $d F_Y (y)$ is absolutely continuous with respect to the Lebesgue measure $dy$, which we assume hereafter to be the case, there exists a probability density function $\rho:\mathcal{D}\to \mathbb{R}$ such that $\rho(y)dy = d F_Y (y)$. Note that the new measure space $(\mathcal{D}, \mathcal{B}(\mathcal{D}), \rho(y)dy)$ is isometric to $(A, \mathcal{F} , P )$ under the random variable $Y$.
\par We define the probability Hilbert space $L^2(A):=\{v:A\to \mathbb{R}:\int_A v^2(\omega)dP(\omega) <\infty \}$ and $L^2_\rho(\mathcal{D}):=\{u :\mathcal{D}\to \mathbb{R} | \int_\mathcal{D} u^2(y) \rho(y) dy <\infty \}$, equipped with the equivalent norms (by noting that $v(\omega)=u(y(\omega))$)
\begin{equation}\label{stochastic_norm}
||v||_{L^2(A)}:=\left( \int_A v^2(\omega) dP(\omega)\right)^{1/2} =\left( \int_\mathcal{D} u^2(y) \rho(y) dy \right)^{1/2}=:||u ||_{L_\rho^2(\mathcal{D})}.  
\end{equation}
\par Let $v:\Omega\times A \to \mathbb{R}$ be a real-valued\emph{ random field}, which is a real-valued random variable defined on $A$ for each $x\in \Omega$. We define the Hilbert space $S(\Omega):=L^2(A) \bigotimes H^1(\Omega)$, equipped with the inner product 
\begin{equation}\label{stochastic_hilbert_inner_product}
(u,v)=\int_A \int_\Omega (uv + \nabla u \cdot \nabla v) \ dx \ dP(\omega)\quad \forall u,v\in S(\Omega),
\end{equation}
where $\nabla$ is the spatial gradient in $\Omega$. The associated norm is defined as $||v||_{S(\Omega)}=\sqrt{(v,v)}.$ 
\par Now we can introduce \emph{stochastic partial differential equations}. Given random vector field $\bmu : A\to \mathbb{R}^p$, our stochastic advection-diffusion problem will be finding a random field $u(x; \bmu(\omega))$ such that
\begin{equation}\label{stochastic_adv_diffusion}
-\varepsilon (\bmu(\omega)) \Delta u(\bmu(\omega)) + \bbe (\bmu(\omega)) \cdot \nabla u(\bmu(\omega)) = f(\bmu(\omega))  \qquad\text{in }\Omega (\bmu(\omega)),
\end{equation}
accompanied by suitable boundary conditions.
%\par As in previous chapter, we have to consider the weak formulation of this problem \eqref{eq:original_problem}, then a discretized version of this weak formulation with FE \eqref{problem_FE}, and we have to remember the RB approach and the greedy algorithm of section \ref{RB}.
%\par \textcolor{red}{FB: rileggendo quello che segue penso ci sia da usare la norma V al quadrato in $\mathbb{E}[...]$ anziche il quadrato della differenza}
\par Now, we want to develop an algorithm that gives more importance to parameters with higher probability of being chosen. The basic idea is to assign different weights to every values of parameter $\bmu\in \mathcal{D}\subset \mathbb{R}^p$ according to a prescribed weight function $w(\bmu) > 0$, and to use them during the procedure of construction of the RB space. The motivation is that when the parameter $\bmu$ has non constant weight function $w(\bmu)$, e.g. stochastic problems with random inputs obeying probability distribution far from uniform type, the weighted approach can considerably attenuate the computational effort for large scale computational problems. The weighted reduced basis method consists of the same elements, namely Greedy algorithm, a posteriori error estimate and \emph{Offline--Online} decomposition, as presented in section \ref{RB}. In this section, we only highlight the new weighted steps.
\par Let $X^\mathcal{N}$ be a high-fidelity approximation space of $X$, equipped with the norm $|||.|||_{\bmu}$ defined in section \ref{sec:rbdet}. Moreover, let us define an equivalent weighted norm
\begin{equation}\label{stochastic_norm_weighted_space}
||u(\bmu)||_w=w(\bmu)||u(\bmu)||_{\bmu}\quad \forall u \in X^\mathcal{N},\forall \bmu \in \mathcal{D},
\end{equation} 
where $w:\mathcal{D}\to \mathbb{R}^+$ is a weighted function taking positive real values, which we assume to be continuous and bounded. We will denote by $X_w$ the space $X$ endowed with $||\cdot||_w$.
\par The Greedy algorithm is thus modified to take the weighting into account, that is to solve an optimization problem in $L^\infty (\mathcal{D};X_w)$: at each step we are seeking a new parameter $\bmu^N\in\mathcal{D}$ such that
\begin{equation}
\bmu^N=\arg\sup_{\bmu\in\Xi_{train}} ||u^\mathcal{N}(\bmu)-u_N(\bmu)||_w,
\end{equation}
where again $u_N$ is the reduced basis approximation of the \emph{truth }solution $u^\mathcal{N}$. Here, $\Xi_{train}$ is the discretized version of the parameter space $\mathcal{D}$. Instead of performing the true error, we use a weighted \emph{a posteriori} error estimator $\Delta_N^w$ such that
\begin{equation}\label{stochastic_Delta_N}
||u^\mathcal{N}(\bmu)-u_N(\bmu)||_w\leq \Delta_N^w(\bmu).
\end{equation}
The choice of the weight function $w(\bmu)$ is aimed by the desire of minimizing the squared norm error of the RB approximation in the space $L^\infty (\mathcal{D};X_w) $, i.e.
\begin{equation}\label{square_RB_w_error}
\begin{split}
\mathbb{E}[||u^\mathcal{N}-u_N||^2]=&\int_A \int_\Omega ||u^\mathcal{N}\bmu(\omega))-u_N(\bmu(\omega))||^2_{\bmu} dx\ dP(\omega)=\\
=&\int_\mathcal{D} \int_\Omega ||u^\mathcal{N}(\bmu)-u_N(\bmu)||^2_{\bmu} \rho(\bmu)\ dx\ d\bmu,
\end{split}
\end{equation}
that we can bound with
\begin{equation}\label{RB_w_error_bounding}
\mathbb{E}\left[||u^\mathcal{N}-u_N||^2\right]\leq \int_\mathcal{D} \Delta_N(\bmu)^2 \rho(\bmu) d\bmu,
\end{equation}
where $\Delta_N$ is the RB error estimator introduced in section \ref{RB}. This motivates us in the choice $w(\bmu)=\sqrt{\rho(\bmu)}$. Finally, we set $\Delta_N^w(\bmu):=\Delta_N(\bmu)\sqrt{\rho(\bmu)}$ \cite{peng_art}.
\par
Another important aspect in the RB algorithm is the choice of the training set $\Xi_{train}$. While in the deterministic case we used Uniform Monte Carlo sampling methods to choose elements from $\mathcal{D}$, in the stochastic context we can use a Monte Carlo sampling according to the distribution $\rho(\bmu)$. We will see in numerical test that this choice is important to improve the convergence of the error.\\
We refer to \cite{peng, peng_art, peng_A} for further details on weighted reduced basis methods.

\subsection{Stabilized weighted reduced basis methods}
In this section we study a variant of the weighted reduced basis method suited for stochastic advection--diffusion equations with high P\'eclet number. In order to do so, we combine the stabilization of advective terms, introduced in section \ref{stab_RB}, to the weighting procedure of section \ref{sec:weighted_RB}.

As in section \ref{section stabilized}, for the moment, we need to add SUPG stabilization terms to the weak form of the problem. This results in the following formulation:
\begin{equation}\label{graetz_SUPG_weighted}
\begin{split}
&\text{find }u^\mathcal{N}(\bmu (\omega)) \in X^\mathcal{N} \text{ s.t. }\\
&a_{stab}(u^\mathcal{N}(\bmu(\omega)), v^\mathcal{N}; \bmu(\omega) ) = F_{stab}(v^\mathcal{N};\bmu (\omega)) \quad v^\mathcal{N}\in X^\mathcal{N} ,\,\forall \omega \in A,
\end{split} 
\end{equation}
where $a_{stab}$ and $F_{stab}$ are defined in section \ref{section stabilized}. The most relevant difference with respect to the previous section is that $\bmu: A\to \mathcal{D}$ is a random vector, instead of being a deterministic parameter.

We test the proposed method with stochastic versions of the previous test cases (\PGP problem \ref{Graetz} and \PFS problem \ref{section:square}). In order to do so, we need to prescribe the distribution of $\bmu$; this will be done for each test case in the following sections. For the sake of exposition results are presented only for the Offline-Online stabilization.
\begin{comment}
In what follows, we will choose the following distribution:
\begin{align}\label{Beta_distr_graetz}
&\mu_1\sim \mu_{1,min} + (\mu_{1,max}-\mu_{1,min})X_1,&X_1\sim \text{Beta}(\alpha_1,\beta_1),\\
&\mu_2\sim \mu_{2,min} + (\mu_{2,max}-\mu_{2,min})X_2,&X_2\sim \text{Beta}(\alpha_2,\beta_2).
\end{align} 
We choose the Beta distribution because it takes values in a compact set, and because we can give more importance to a certain subset of the range (for example, the one with higher P\'eclet number). The prescribed values of $\alpha_1,\beta_1, \alpha_2,\beta_2$ will be clarified in each test case.
\end{comment}
%In figure \ref{esempio_beta} we can see the density function of a $\text{Beta}(4,2)$ distribution.
%\begin{figure}[h]
%\centering
%\includegraphics[trim=0cm 0cm 0cm 0cm,clip=true, scale=0.45]{Immagini/Stochastic/Beta4_2.pdf} 
%\caption{Beta(4,2)} \label{esempio_beta}
%\end{figure}

\subsubsection{Numerical test: Poiseuille--Graetz problem}\label{stochastic_graetz}
For \PGP problem, we consider the range $\mathcal{D}=\left[ 10^1,10^6\right]\times \left[0.5,4 \right]$ for the parameter $\bmu$. To give more importance to parameter with $\mu_1\approx 10^{5}$, we use $X_1\sim \text{Beta}(4,2)$ and $\mu_1 \sim 10^{1+5\cdot X_1}$, while $X_2\sim \text{Beta}(3,4)$ and $\mu_2\sim 0.5 + 3.5 X_2$. We choose the Beta distribution because it takes values in a compact set\footnote{The weighted approach would work as well for an unbounded (e.g. Gaussian) distribution. We use a Beta distribution in order to be able to present the comparison between a weighted and the classical approach. The latter would not be possible for Gaussian random variables, unless the parameter domain is cut. Such cut would be somehow arbitrary, since the classical approach does not exploit the underlying probability distribution.}, resulting in $(\mu_1, \mu_2) \in \mathcal{D}$.
%In previous chapters we, implicitly used this algorithm with a uniform distribution over $\mathcal{D}$.\\

\begin{figure}[h]
\centering
\subfigure[Error Comparison]{\includegraphics[trim=0cm 0cm 1.5cm 0cm,clip=true, scale=0.38]{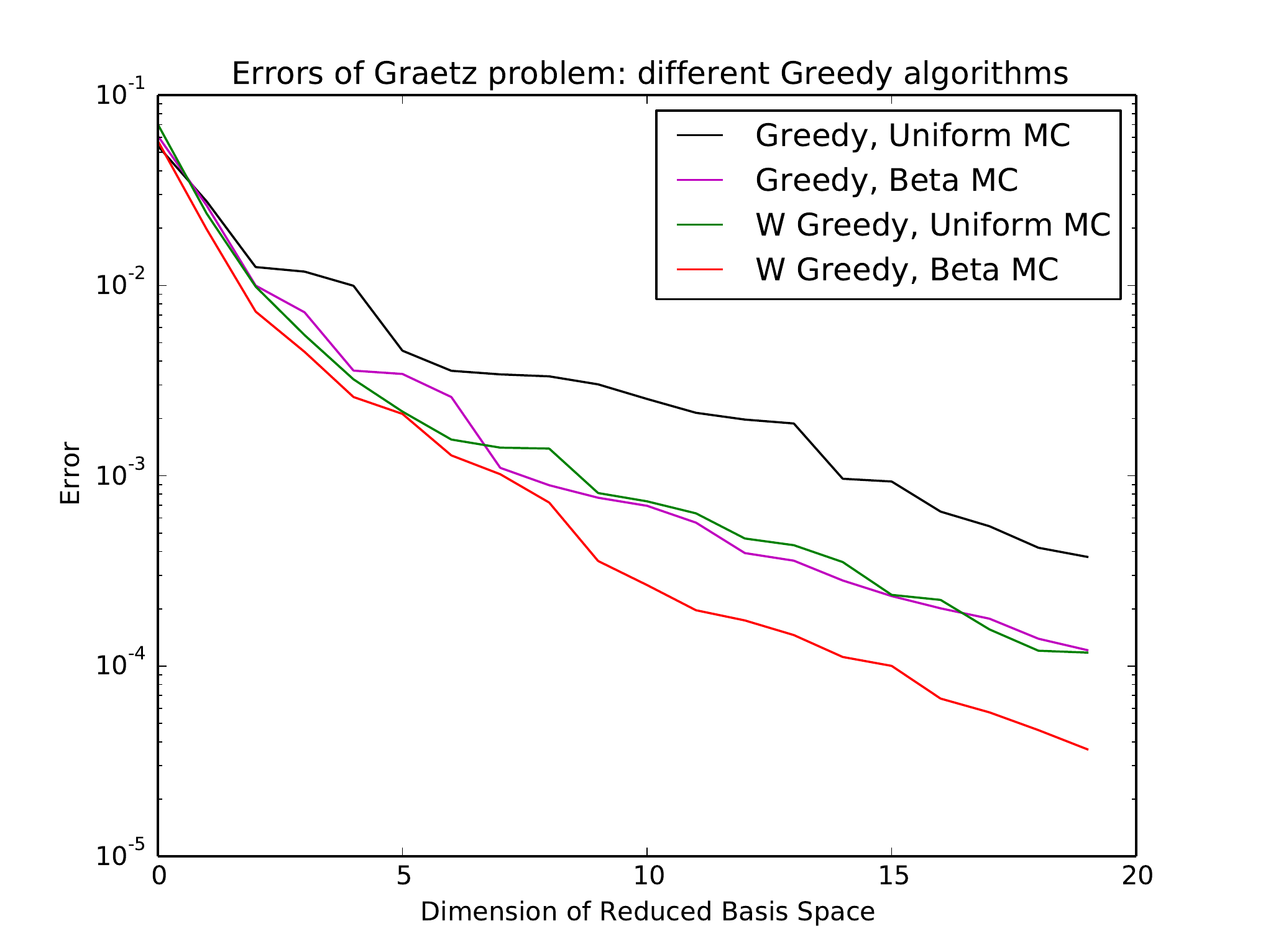} }
\subfigure[$\Delta_N$ comparison]{\includegraphics[trim=0cm 0cm 1.5cm 0cm,clip=true, scale=0.38]{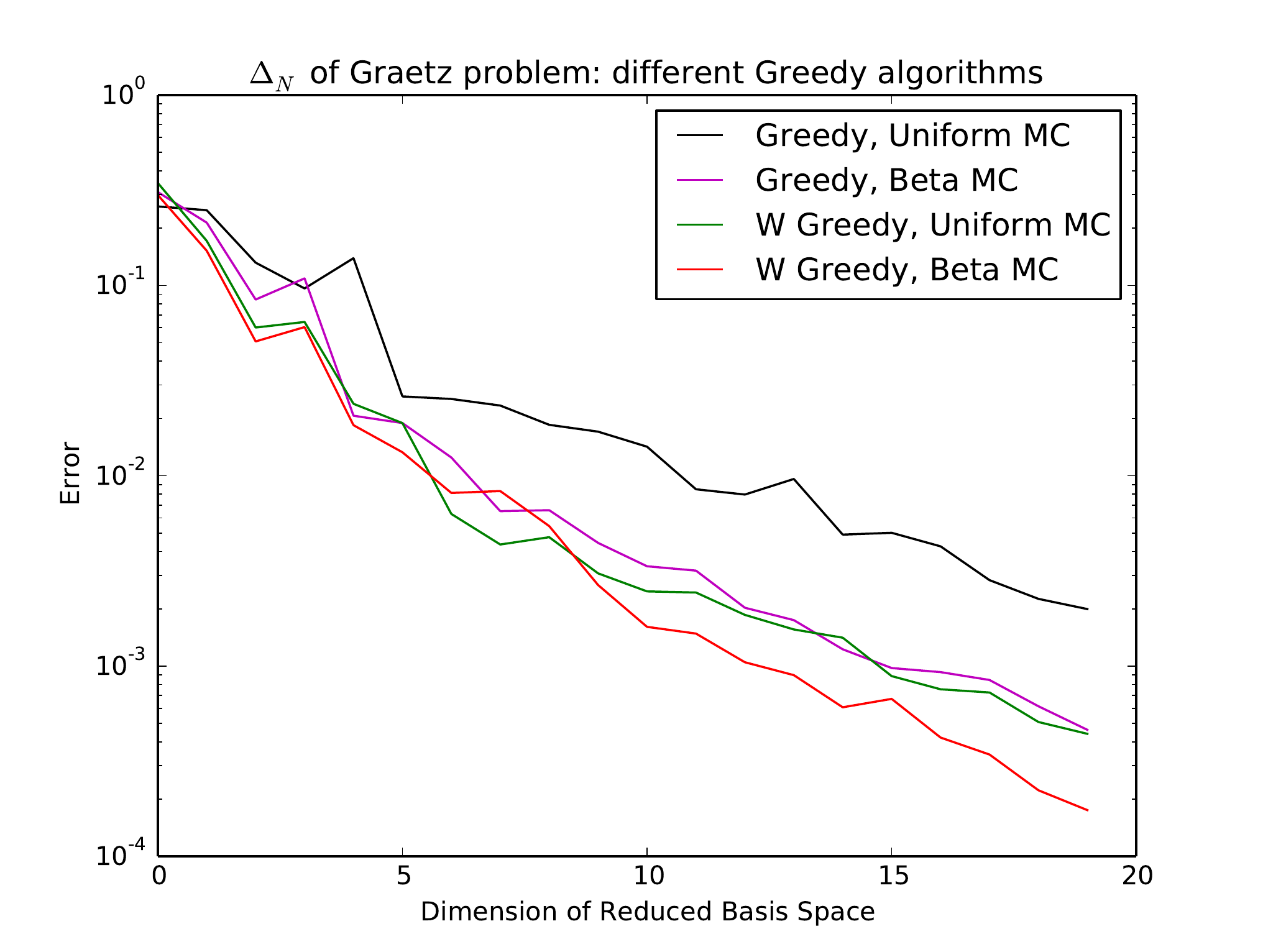} }
\caption{Greedy algorithms comparison for Graetz problem} \label{stochastic_graetz_comparison_greedy}
\end{figure}

We compare next the performance of the reduction method for the different choices that we have discussed in section \ref{stochastic}, namely related to using weighted or standard Greedy algorithm, and the sampling of the training set $\Xi_{train}$. We present in figure \ref{stochastic_graetz_comparison_greedy} numerical results for four different cases: 
\begin{enumerate}
\item Classical Greedy with Uniform Monte Carlo sampling (black line);
\item Classical Greedy with Beta Monte Carlo sampling (purple line);
\item Weighted Greedy with Uniform Monte Carlo sampling (green line);
\item Weighted Greedy with Beta Monte Carlo sampling (red line).
\end{enumerate}

We used 200 samples for $\Xi_{train}$ in each algorithm during the offline stage. We can see in figure \ref{stochastic_graetz_comparison_greedy} the comparison between the average errors and the average $\Delta_N$ between these algorithms for a test set of size 100, with the same distribution as the training set. The results show that both weighting and a correct sampling are essential to obtain the best convergence results \cite{tesi_luca,VenturiBallarinRozza2018}.
Indeed, putting together these two aspects we get the best results, reaching an error that is one tenth of the error of the classical Greedy algorithm on uniform distribution.
\par
In a similar way, instead of computing the average of the errors on the test set, we can also compute the mean of the error in a probability sense, i.e.
\begin{align}
\mathbb{E}[|||u^\mathcal{N}(\bmu)-u_N(\bmu)|||_{\bmu}] & = \int_\mathcal{A} |||u^\mathcal{N}(\bmu(\omega))-u_N(\bmu(\omega))|||_{\bmu(\omega)} dP(\omega)\label{error_stochastic}\\
& = \int_\mathcal{D} |||u^\mathcal{N}(\bmu)-u_N(\bmu)|||_{\bmu} \rho(\bmu) d\bmu,\label{error_stochastic_2}
\end{align}
that we can approximate using some quadrature method. In particular, we will use Monte Carlo method, i.e. we approximate \eqref{error_stochastic} with
\begin{equation}\label{error_montecarlo}
\mathbb{E}[|||u^\mathcal{N}(\bmu)-u_N(\bmu)|||_{\bmu}]\approx \frac{1}{M} \sum_{i=1}^M |||u^\mathcal{N}(\bmu_i)-u_N(\bmu_i)|||_{\bmu_i},
\end{equation}
where $\bmu_i,\,i=1,\dots ,M$ are random parameters in the testing test drawn from a Beta distribution,
while we approximate \eqref{error_stochastic_2} with
\begin{equation}\label{error_montecarlo_2}
\mathbb{E}[|||u^\mathcal{N}(\bmu)-u_N(\bmu)|||_{\bmu}]\approx \frac{1}{M} \sum_{j=1}^M |||u^\mathcal{N}(\bmu_j)-u_N(\bmu_i)|||_{\bmu_j} \rho(\bmu_j),
\end{equation}
where $\bmu_j,\,i=1,\dots ,M$ are drawn from a Uniform distribution (on the same support) instead.

Results are nevertheless similar to the ones presented in figure \ref{stochastic_graetz_comparison_greedy}, and the same conclusions can be drawn. For instance, the probabilistic mean of the errors in the classical Greedy method with uniform sampling and the weighted reduced one with Beta sampling are $4.5485\cdot 10^{-4}$ and $1.2807 \cdot 10^{-4}$, respectively. 

\subsubsection{Numerical test: propagating front in a square}\label{stochastic_square}
We can proceed in the same way for the \PFS problem of section section \ref{section:square}. In this section, the parameter range $\mathcal{D}$ is $\left[10^4, 10^5 \right] \times \left[0, 1.5\right]$. Also in this case $\mu_1$ and $\mu_2$ depend on randomly distributed Beta variables, i.e. $\mu_1\sim 10^4 + 9\cdot 10^4\cdot X_1$ and $\mu_2\sim 1.5\cdot X_2$, where $X_1\sim \text{Beta}(3,4)$ while $X_2\sim \text{Beta} (4,2)$.

\begin{figure}[h]
\centering
\includegraphics[trim=0cm 0cm 1.5cm 0cm,clip=true, scale=0.35]{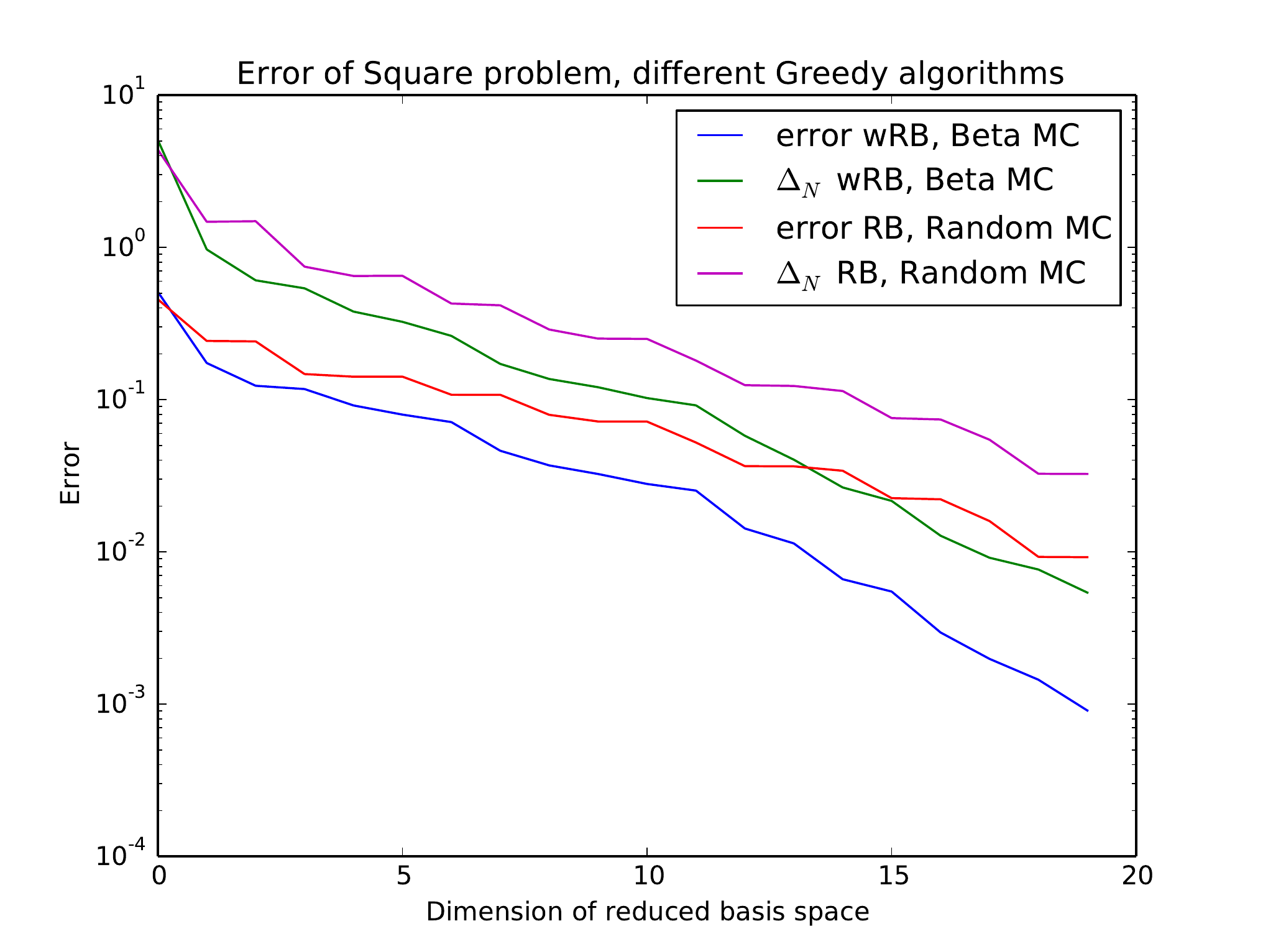}
\caption{Greedy algorithms comparison for \PFS problem} \label{stochastic_square_comparison_greedy}
\end{figure}

As for the previous test case we compare the classical Greedy method with Uniform Monte Carlo to the weighted reduced basis method with Beta Monte Carlo distribution. The comparison, shown in figure \ref{stochastic_square_comparison_greedy}, provides results which are very similar to \PGP problem. Indeed, the weighted RB method with Beta distribution is converging faster than the classical one. Also the mean errors in the probabilistic sense of \eqref{error_montecarlo} show a similar behavior: for a reduced basis space of dimension $N=20$, the stabilized weighted method with Beta distribution produces a mean error of $1.7803 \cdot 10^{-3}$, while the classical approach gives a mean error of $7.9362 \cdot 10^{-3}$.

\subsection{Selective online stabilization of weighted reduced basis approach}
In this section we want to optimize computational costs in the \emph{Online} phase of RB method. Indeed, stabilization procedure can lead to an increase in the number $Q_a$ and/or $Q_f$ of affine terms, which in turn may lead to larger online times required for the assembly of the linear system or for the evaluation of the error estimator. In this section we propose a procedure to selectively enable online stabilization only when required. In the whole section we keep the reduced basis produced in the previous section for $N=20$.

%We know that errors computed with the \emph{Offline--only} stabilized method are often very large, compared with \emph{Offline--Online }stabilization (for example for Graetz problem in figure \ref{graetz_error_comp_onl_off}), especially for high P\'eclet numbers. Working in context of stochastic partial differential equation and using the error \eqref{error_stochastic}, we want to combine the \emph{Offline--only} and the \emph{Offline--Online }approach to reduce computational costs and not to worsen the error too much.

\subsubsection{Numerical test: Poiseuille--Graetz problem}
Let us consider first the \PGP example, with Beta distribution over parameter $\bmu$, similarly to section \ref{stochastic_graetz}. 
In what follows, we assume that $\mu_1 \in [10, 10^6]$,  $\mu_1 \sim 10^{1+5\cdot X_1}$ where $X_1\sim \text{Beta}(5, 3)$. To simplify the discussion of the results we further assume that $\mu_2 \equiv 1$.

\begin{figure}[h]
\centering
\subfigure[Beta(5,3) distribution]{\includegraphics[trim=0cm 0cm 0cm 0.5cm,clip=true, scale=0.35]{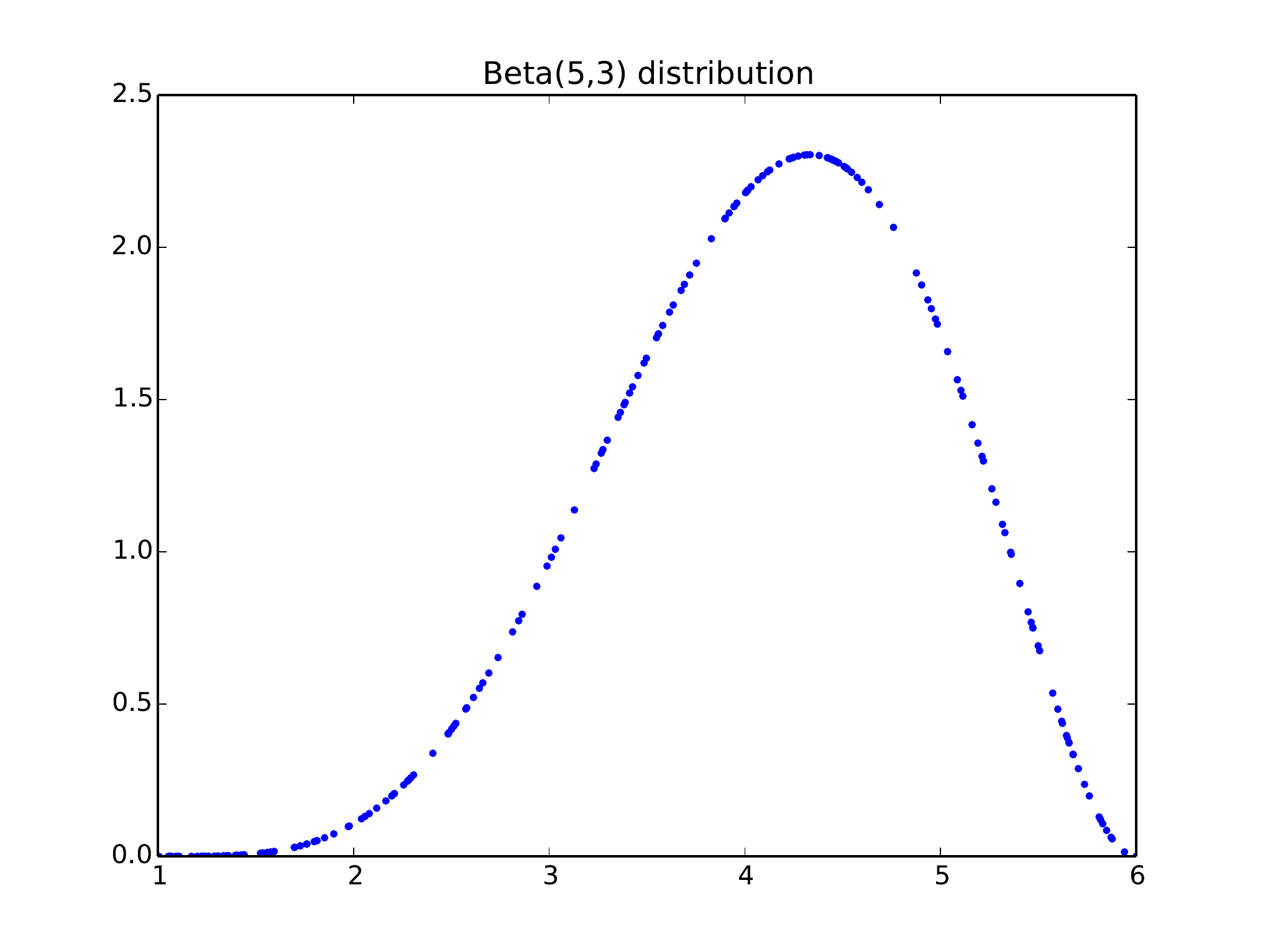} \label{stochastic_beta_5_3} }
\subfigure[Errors with stabilization Offline and Offline-Online]{\includegraphics[trim=0cm 0cm 0cm 0.5cm,clip=true, scale=0.35]{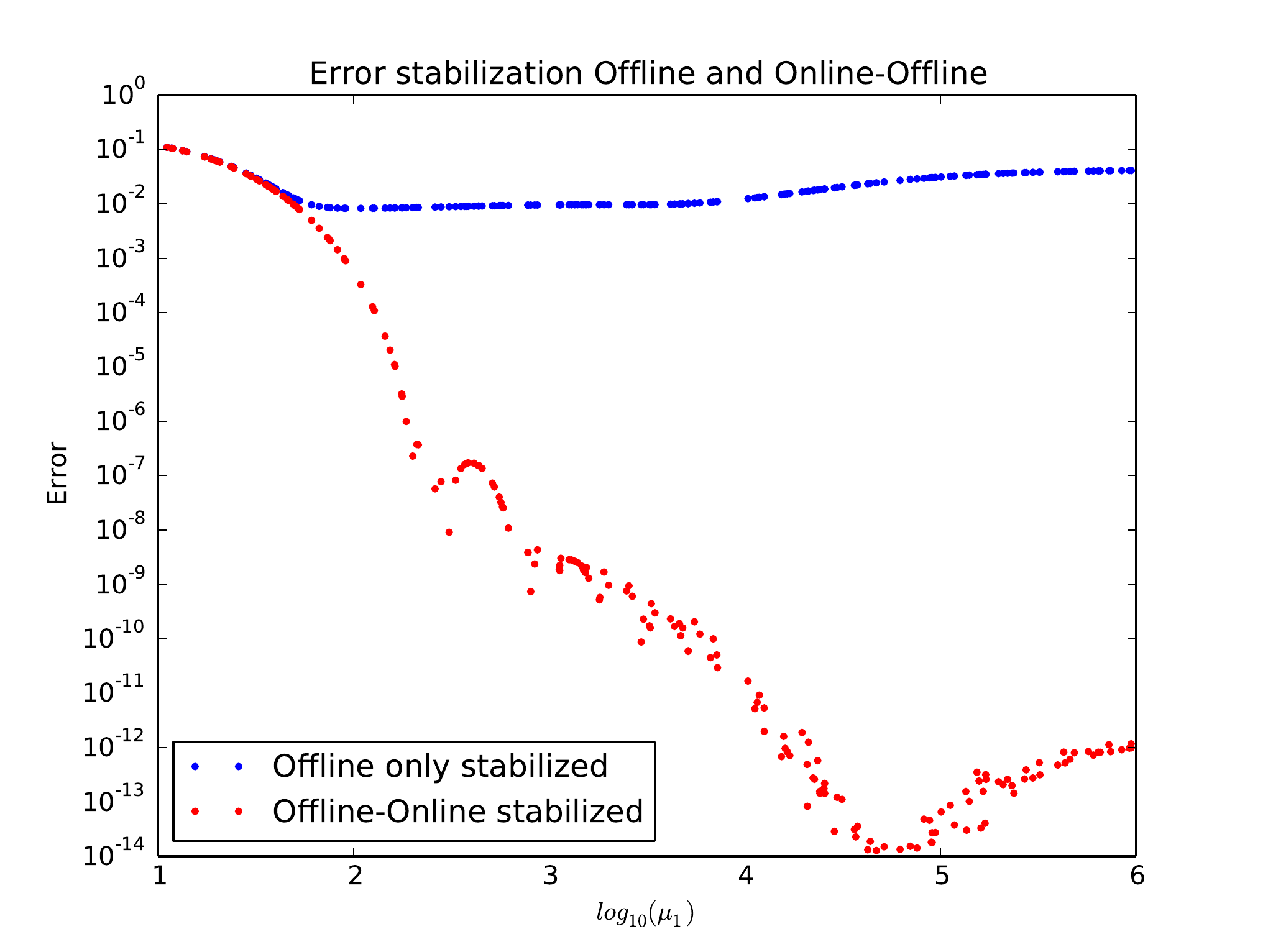} \label{stochastic_stab_off_onloff}}
\caption{Error and density of Uniform Monte Carlo test set}
\end{figure}

%To sum up what we have done in the previous section, we have taken $\bmu$ as a random vector, we have taken the stochastic equation \eqref{stochastic_adv_diffusion}, its weak formulation and the discretized version. Then, we applied the \emph{Offline} stage of the weighted RB algorithm (section \ref{sec:weighted_RB}) to SUPG stabilized version of this equation \eqref{graetz_SUPG_weighted}.\\
While carrying out the online stage of the proposed stabilized weighted reduced basis method, we can choose whether to apply online stabilization or not. 
Figure \ref{stochastic_stab_off_onloff} shows the resulting error on a test set (that we have taken with a Uniform Monte Carlo sampling), sorted by increasing values of $\mu_1$, considering both options. We can observe that for low P\'eclet number ($\mu_1\leq 10^2$), \emph{Offline-Online} stabilization and \emph{Offline only} stabilization produce very similar results. Thus, we would prefer the less expensive \emph{Offline only} stabilization procedure. There the error is high, because the samples selected from the weighted Greedy in the Offline phase are all concentrated where the density of probability is higher (high P\'eclet). For this reason the low P\'eclet number zone is bad represented. Moreover, in the regions where the density of $\bmu$ is very small, even a large error would be less relevant in terms of the probabilistic mean error \eqref{error_stochastic}.
So, we should consider the idea of enabling the more expensive online stabilization only for parameters with high density (which would affect more the mean error) or parameters with large P\'eclet numbers (were the more expensive assembly is fully justified by the convection dominated regime).

\begin{figure}[h]
\centering
\includegraphics[trim=0cm 0.45cm 0cm 0.5cm,clip=true, scale=0.35]{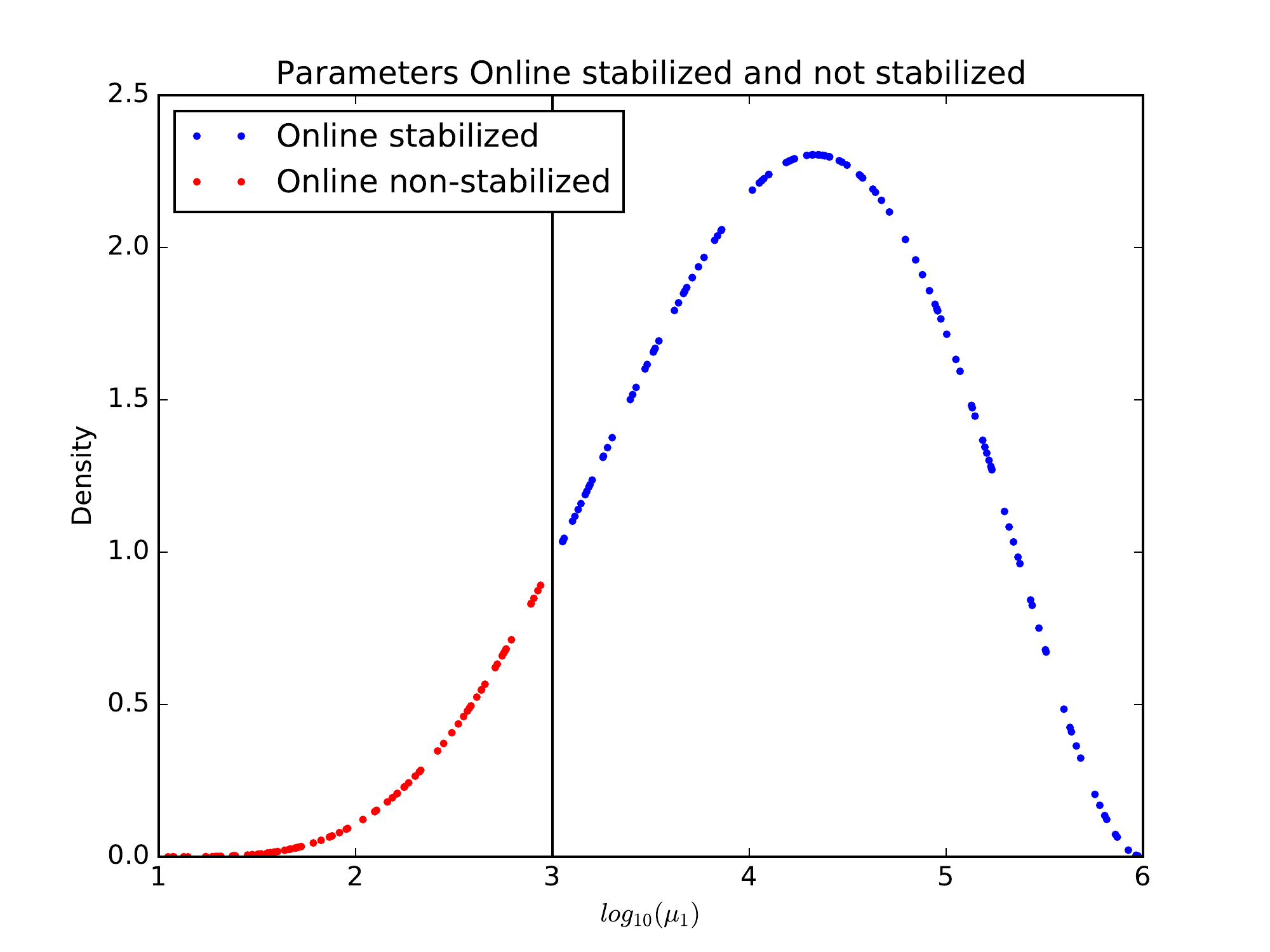} 
\caption{P\'eclet discriminant, black line is the P\'eclet threshold} \label{distribution_peclet_off_onl}
\end{figure}

Let us start considering the case where we want to stabilize \emph{Online} solutions depending on P\'eclet numbers. First, we establish a threshold at a certain P\'eclet number $\widetilde{\mu}_1$. For parameters $\mu_1>\widetilde{\mu}_1$ we will use both \emph{Online} and \emph{Offline} stabilization, while for parameter $\mu_1\leq \widetilde{\mu}_1$ we will use only \emph{Offline} stabilization. See figure \ref{distribution_peclet_off_onl} for a graphical representation for $\widetilde{\mu}_1 = 10^3$.\\
For different thresholds $\widetilde{\mu}_1$ we can compute the error in sense of \eqref{error_stochastic}, as we can see in the following table.
 %After this, we can compute the error in sense of \eqref{error_stochastic} and we can see that, if during the \emph{Online} stage we do not stabilize every solution, we have an error of 0.021128, while stabilizing everything we get an error of $7.967 \cdot 10^{-4}$. 
$$
\begin{array}{c|c|c}
\text{Threshold $\widetilde{\mu}_1$} & \text{Error} & \text{Percentage non-stabilized}\\
\hline
10^{1}     &  7.9673 \cdot 10^{-4} & 0\%  \\
10^{1.5} & 8.0704 \cdot 10^{-4} & 10\%  \\
10^{2} & 10.0060 \cdot 10^{-4} & 20\% \\
10^{2.5} & 18.2806 \cdot 10^{-4} & 33\%\\
10^3 & 33.4593 \cdot 10^{-4} & 45\%\\
10^{6} & 0.021128 &100\%
\end{array}
$$
Considering that the best attainable error was of $7.967\cdot 10^{-4}$, we can say that until $\widetilde{\mu}_1 = 10^{2}$ we are not worsening considerably the error (less than an order of magnitude). At the same time, we can save online time on the assembly of terms related to stabilization coefficient for 20\% of our test set (that was uniformly distributed).

\begin{figure}[h]
\centering
\includegraphics[trim=0cm 0.45cm 0cm 0.5cm,clip=true, scale=0.35]{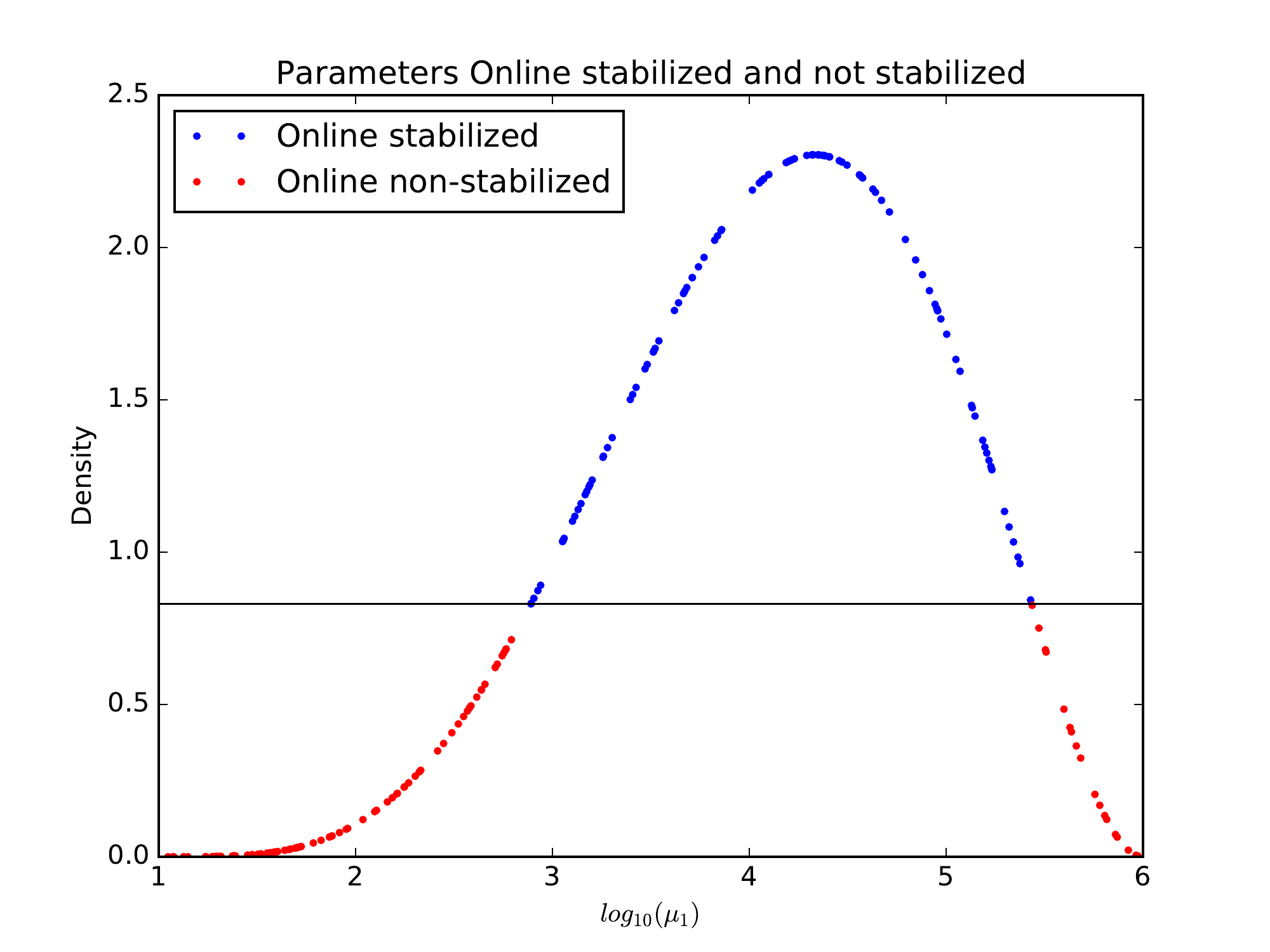} 
\caption{Density discriminant, black line is the density threshold} \label{distribution_density_off_onl}
\end{figure}

The other natural gauge to decide whether to stabilize \emph{Online}, or not, is the density $\rho(\bmu)$. 
Let $\widetilde{\nu}$ be a prescribed tolerance; we will not stabilize parameters $\bmu$ on the tail $I$ of the distribution such that 
\begin{equation}
\int_I \rho(\bmu) d\bmu = \widetilde{\nu},
\end{equation}
where $I$ is a set $\{ \bmu : \rho(\bmu)\leq \widetilde{\rho} \}$ for some suitable $\widetilde{\rho}$ which can be easily found numerically as a function of $\widetilde{\nu}$.
In figure \ref{distribution_density_off_onl} we can see an example for $\widetilde{\nu} = 10\%$.

In the following table, we summarize some results for different thresholds $\widetilde{\nu}$ (and, correspondigly, $\widetilde{\rho}$).
$$
\begin{array}{c|c|c|c}
\text{Threshold }\widetilde{\nu} & \text{Threshold }\widetilde{\rho} & \text{Error} & \text{Percentage non-stabilized}\\
\hline
0     & 0       & 7.9673 \cdot 10^{-4} & 0\%  \\
0.001 & 0.02233 & 9.3222 \cdot 10^{-4} & 15\%  \\
0.002 & 0.04423 & 9.6456 \cdot 10^{-4} & 17\%  \\
0.005 & 0.09094 & 14.7861 \cdot 10^{-4} & 21\%  \\
0.01  & 0.13877 & 15.9482 \cdot 10^{-4} & 25\%  \\
0.02  & 0.21433 & 25.6017 \cdot 10^{-4} & 30\%  \\
0.05  & 0.38244 & 49.1931 \cdot 10^{-4} & 38\%  \\
0.1   & 0.89068 & 66.7488 \cdot 10^{-4} & 45\%  \\
1     & \infty  & 0.021128 & 100\%  
\end{array}
$$
We have that errors computed using density discriminant are less accurate than ones computed with P\'eclet discriminant. Indeed, for the same percentage of non-stabilized solution (for example 45\%) we have bigger errors in density discriminant approach ($66 \cdot 10^{-4}$ instead of $33 \cdot 10^{-4}$). This is due to the enormous difference between \emph{Online} stabilized and \emph{Online} non--stabilized solution for high P\'eclet numbers (figure \ref{stochastic_stab_off_onloff}), with the latter resulting in considerably larger errors.

\subsubsection{Numerical test: propagating front in a square}
Let us now consider the \PFS problem with fixed $\mu_1 \equiv 10^5$, while $\mu_2\sim 0.5 + 3.5 X_2 \in [0, 1.5]$ where $X_2\sim \text{Beta}(4, 2)$. We have decided to fix the P\'eclet number since results in section \ref{section:square} show that the solution is most sensible to the parameter $\mu_2$, which represents the angle of the propagating front. 

\begin{figure}[h]
\centering
\subfigure[]{\includegraphics[trim=0cm 0.45cm 0cm 0.5cm,clip=true, scale=0.35]{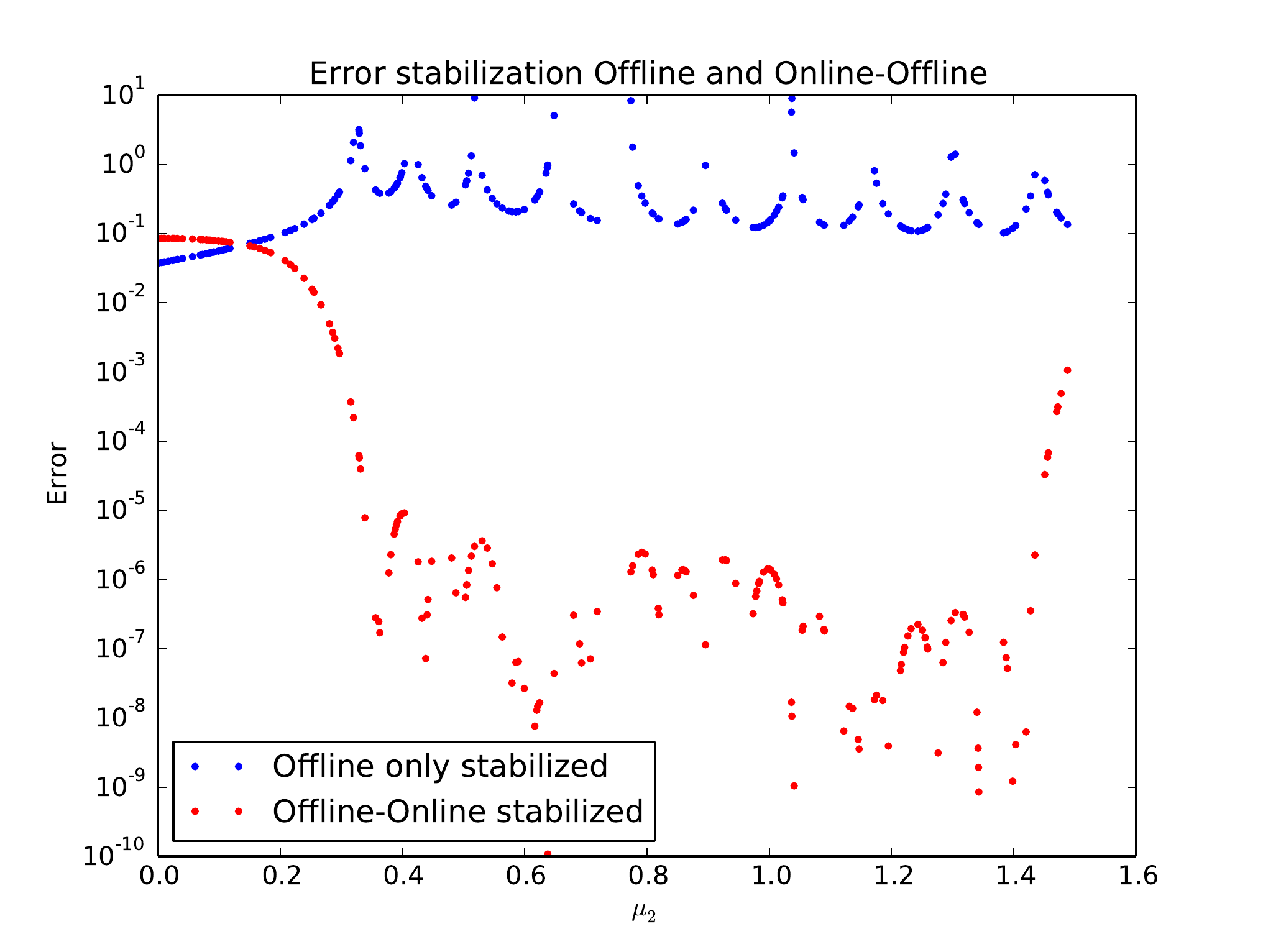} \label{sqaure_error_online_offline}}\\
\subfigure[]{\includegraphics[trim=0cm 0.45cm 0cm 0.5cm,clip=true, scale=0.35]{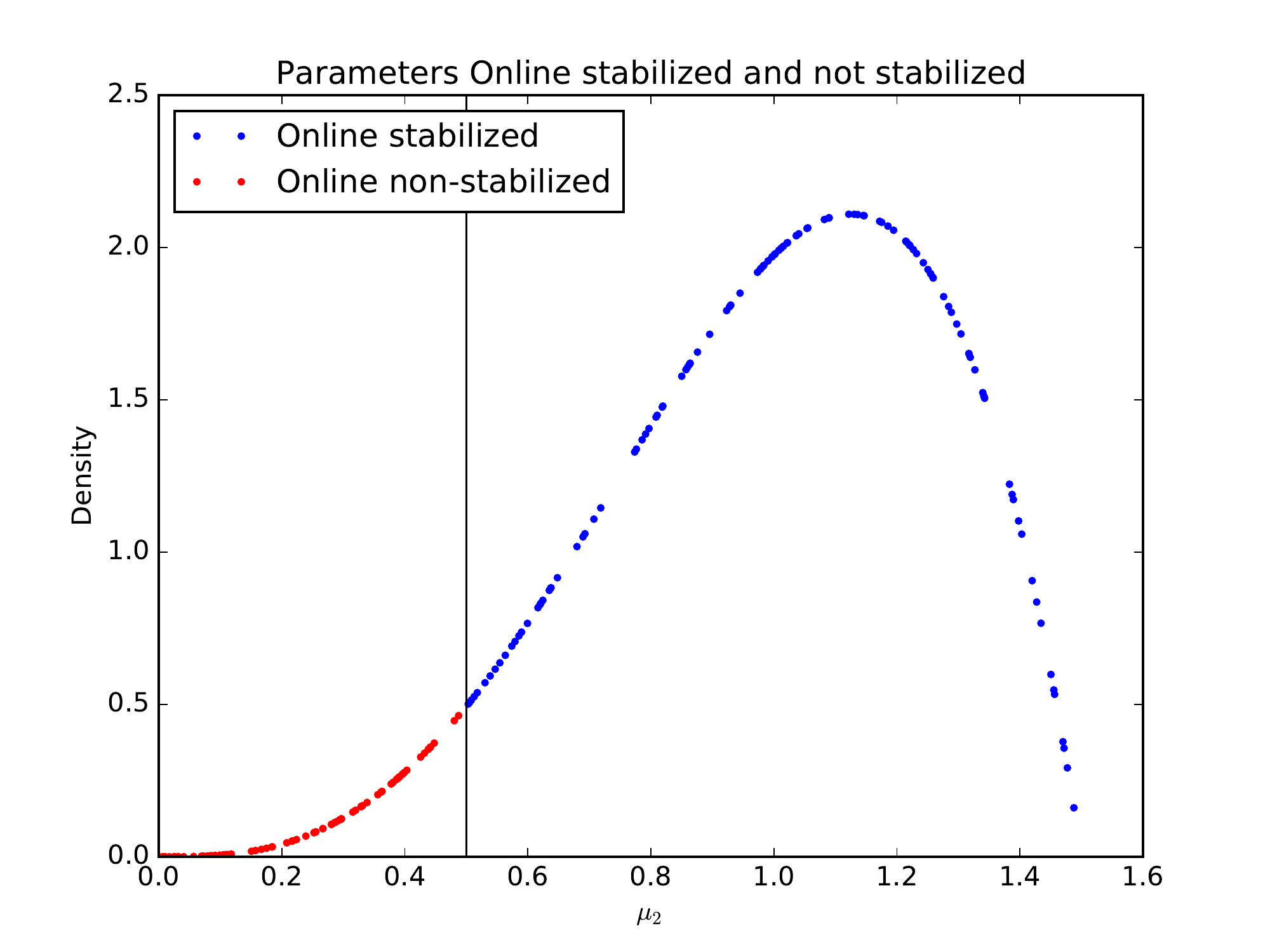} \label{square_threshold_angle}}
\subfigure[]{\includegraphics[trim=0cm 0.45cm 0cm 0.5cm,clip=true, scale=0.35]{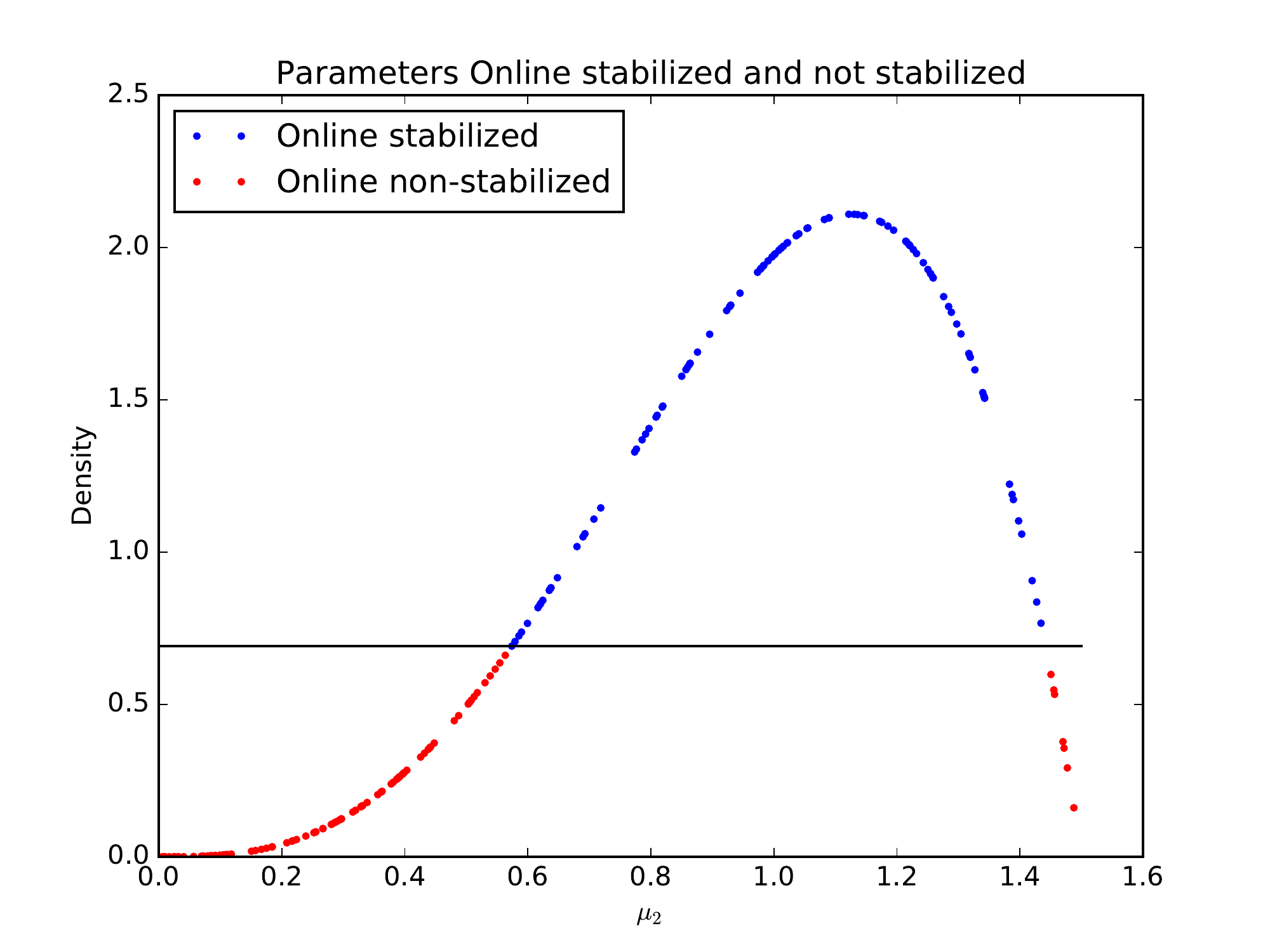}  \label{square_threshold_density}}
\caption{(a) Errors with stabilization Offline and Offline-Online; (b) angle discriminant, black line is the angle threshold; (c) density discriminant, black line is the density threshold }
\end{figure}

Errors for \emph{Online} stabilized and \emph{Online} not stabilized solutions over a Uniform Monte Carlo test set of 200 elements are provided in figure \ref{sqaure_error_online_offline}, for increasing values of $\mu_2$.
We can notice that \emph{Offline--Online} stabilized errors of solutions with small angles (figure \ref{sqaure_error_online_offline}, $\mu_2 \lesssim 0.2$) are bigger than \emph{Offline--only} stabilized errors. This is due to the fact that the density of that region of the parameter range is very small and thus the weighted Greedy algorithm picks very few parameters in that region. In a similar way, we also notice that solutions for $\mu_2 \approx 1.5$ are not well approximated. Indeed, in the \emph{Offline only} stabilized case the lack of stabilization is badly affecting the reduced order solution for any $\mu_2 \gtrsim 0.2$, while in the \emph{Offline-Online} stabilized case the low density of $\mu_2 \gtrsim 1.4$ leads the weighted reduced basis selection to choose few parameters $\mu_2 \approx 1.5$ during the offline stage.

Thus, in a similar way to the previous test case, we propose selective online stabilization criteria, either depending on a threshold on the parameter (the angle $\mu_2$ in this case, rather than the P\'eclet number) or on the probability distribution. Let us start from a discussion of the former choice, leading to online stabilize for angles greater than a certain threshold $\widetilde{\mu}_2$ (see e.g.\ figure \ref{square_threshold_angle}).
%The error in sense of \eqref{error_stochastic} of all stabilized online solutions is of 0.01416, while with all not stabilized is of 0.82998. 
The error for different thresholds $\widetilde{\mu}_2$ is tabulated as follows:
$$
\begin{array}{c|c|c}
\text{Threshold } \widetilde{\mu}_2 & \text{Error} & \text{Percentage non-stabilized}\\
\hline
0     & 0.01416  & 0\%  \\
0.1 & 0.01400 & 6\% \\
0.2 & 0.01506 & 16\%\\
0.3 & 0.04056 & 23\%\\
0.4 & 0.11810 & 30\%\\
0.5 & 0.20365 & 37\%\\
1.5 & 0.82998 &100\%
\end{array}
$$
We can observe that at the beginning the error is decreasing as the threshold increases, while it slowly increases after a critical angle between $0.1$ and $0.2$. Due to this, we consider a threshold $\widetilde{\mu}_2 = 0.2$ to be optimal in order not to increase the error and save 16\% of online stabilization computations.

As for \PGP example, we can also test a criterion based on a density threshold (see e.g.\ figure \ref{square_threshold_density}). In the following table, we are showing different errors for different density thresholds.

$$
\begin{array}{c|c|c|c}
\text{Threshold }\widetilde{\nu}& \text{Threshold }\widetilde{\rho} & \text{Error} & \text{Percentage non-stabilized}\\
\hline
0     & 0       & 0.01416 & 0\%  \\
0.001 & 0.02271 & 0.01400 & 13\%  \\
0.002 & 0.04600 & 0.01506 & 16\%  \\
0.005 & 0.10237 & 0.02269 & 20\%  \\
0.01  & 0.13598 & 0.04658  & 25\%  \\
0.02  & 0.26309 & 0.11158  & 30\%  \\
0.05  & 0.51855 & 0.20613 & 38\%  \\
0.1   & 0.72557 & 0.32034 & 46\%  \\
1     & \infty  & 0.82998 & 100\%  
\end{array}
$$
In this case, a negligible increase of the error is obtained for $\widetilde{\nu} = 0.002$, allowing to save more than 15\% of stabilized \emph{Online} computations. Further computational savings can be obtained for $\widetilde{\nu} = 0.01$, up to 25\%, at the expense of a larger error. We notice that in this case both criteria give similar results: this is due to the fact that errors are large for both \emph{Offline only} and \emph{Offline-Online} stabilization methods when $\mu_2$ is large or where density $\rho$ is small.

\begin{remark}
Let $I$ be the region of the parameter space where \emph{Offline only} stabilized solution is selected, and let $\mathcal{D} \setminus I$ denote the complement region in which the \emph{Offline-Online} stabilized method is queried. Let $u_{N}^{I}(\bmu)$ denote the corresponding reduced order solution for $\bmu \in I$, and similarly $u_{N}^{\mathcal{D} \setminus I}(\bmu)$ for $\bmu \in \mathcal{D} \setminus I$. To ease the notation, we will denote the online solution by $u_{N}(\bmu)$ when no confusion arises. 

The selective procedure for online stabilization can be automatically tuned according to a prescribed tolerance on the probabilistic mean error $\mathbb{E}\left[|||u^\mathcal{N}(\bmu)-u_N(\bmu)|||_{\bmu}\right]$. In order to estimate the mean error, we recall the standard error estimation \eqref{def:error_estimator} for $\bmu \in \mathcal{D} \setminus I$, and the following error estimation
\begin{equation}
\begin{split}
|||u_N^I(\bmu)-u^\mathcal{N}(\bmu)|||_{\bmu} \leq & \Delta_{N}^{I}(\bmu):=h_{max}(\bmu) C(\bmu) ||\bbe \cdot \nabla u^{\mathcal{N}}(\bmu)||_{L^2(\Omega_p(\bmu))}+\\
& +(1+h_{max}(\bmu) C(\bmu)^2 ||\bbe ||_{L^\infty (\Omega_p (\bmu))})\varepsilon^*.
\end{split}
\end{equation}
for $\bmu \in I$ \cite{pacciarini_a}, where $C(\bmu)$ is the constant of the equivalence between $H^1$ and $|||\cdot|||_{\bmu}$ norms, $h_{max}$ is the maximum mesh size, while $\varepsilon^*$ is the tolerance of the Greedy algorithm \cite{pacciarini_a}.

Thus, combining these two error estimators, we get that
\begin{equation}
\mathbb{E}\left[|||u^\mathcal{N}(\bmu)-u_N|||_{\bmu}\right] \leq (1-\widetilde{\nu} )\max_{\bmu \in \mathcal{D}\setminus  I} \Delta_N(\bmu)+\widetilde{\nu} \max_{\bmu \in I} \Delta_N^{I}(\bmu).
\end{equation}
which, for a given tolerance $\widetilde{e}$ on the mean error, allows us to compute $\widetilde{\nu}$ such that 
$$(1-\widetilde{\nu} )\max_{\bmu \in \mathcal{D}\setminus  I} \Delta_N(\bmu)+\widetilde{\nu} \max_{\bmu \in I} \Delta_N^{I}(\bmu) < \widetilde{e}.$$
\end{remark}
\begin{remark}
We remark that this selective approach for online stabilization is peculiar of stochastic problems. Indeed, it is the density distribution and the 
relative importance of each sample in the computation of the probabilistic mean that drives the selection process. Such a weighting is lacking
in a deterministic setting, being all samples equally probable during the online stage.
\end{remark}

\section{Stabilized weighted reduced basis method for stochastic parabolic equations}\label{parabolic_RB}
In this section we extend our investigation to stochastic time dependent advection--diffusion equations. 
%The RB method for time dependent problems has been already studied in several works, e.g. \cite{gelsomino, reduced_basis, manzoni_45, rozza_nguyen}.
Stabilization of advection diffusion parabolic equations with high P\'eclet number have been studied in several works with different stabilization methods \cite{brooks}. We will adapt SUPG stabilization for FE methods on parabolic equations to RB method, as suggested in \cite{tesi_paolo, pacciarini_a, pacciarini_c, pacciarini_b}. The reduction will employ a POD-Greedy procedure \cite{Haasdonk,Calcolo,manzoni_45} during the offline stage. We refer to \cite{spannring2018phd,spannring2017weighted} for very recent weighted RB variants for stochastic heat equations.
%\par As the previous section, this one is a prerequisite to section \ref{stochastic}, where we will deal with stochastic equations with random input parameters, defined by a prescribed random variable.
%\subsection{Reduced basis methods for linear parabolic equations}\label{parabolic_RB}

Like for stochastic elliptic equations, we define a \emph{parameter domain} $\mathcal{D}$ as a closed subset of $\mathbb{R}^p$ and we call $\bmu$ a random field with values in $\mathcal{D}$. Again, let $\Omega$ be a bounded open subset of $\mathbb{R}^d\,(d=1,2,3)$ with regular boundary $\partial \Omega$ and let $X$ be a functional space such that $H^1_0(\Omega) \subset X\subset H^1(\Omega).$ For each outcome $\omega \in A$, and corresponding realization $\bmu(\omega) \in \mathcal{D}$, we define the continuous, coercive bilinear form $a$ and the continuous, bilinear, symmetric form $m$ such that satisfy the \emph{affinity} assumption like \eqref{affine_a} and a linear form $F$ which satisfies the \emph{affine} assumption \eqref{affine_F}. Let us finally denote the time domain as $I=\left[ 0, T \right]$, where $T$ is the final time.

We can now define the weak form of the continuous stochastic problem:
\begin{equation}\label{parabolic_problem_0}
\begin{split}
&\text{find }u(t;\bmu(\omega)) \in X, \quad \forall t \in I, \quad \forall \omega \in A,\quad \text{continuous in }t \text{ s.t.}\\
&m(\partial_t u(t;\bmu(\omega)),v)+a(u(t;\bmu(\omega)),v;\bmu(\omega))=g(t)F(v;\bmu(\omega))\quad \forall v \in X,\quad \forall t \in I, \quad \forall \omega \in A\\
&\text{given the initial value } u(0;\bmu(\omega))=u_0\in L^2(\Omega)
\end{split}
\end{equation}
where $g:I\to \mathbb{R}$ is a \emph{control function} such that $g\in L^2(I)$. We choose a right hand side of the form $g(t)F(v;\bmu)$, as usual in the RB framework \cite{grepl, manzoni_45}, in order to ease the \emph{Offline--Online }computational decoupling. 
%In particular, the function $g$ can be specified during the \emph{Online }evaluation of the RB solution.
\subsection{Discretization and RB formulation}
To discretize the time--dependent problem \eqref{parabolic_problem_0} we follow the approach used in \cite{grepl, reduced_basis, nguyen, manzoni_45}, that is to use finite differences in time and FE in space discretization \cite{quarteroni_sacco}.
We start by discretizing the spatial part of the problem (resulting in a mesh denoted by $\mathcal{T}_h$) and the temporal part (resulting in discrete time steps $\{t_j=j\cdot \Delta t\}_{j=0}^J$). We thus define the FE \emph{truth} approximation space $X^\mathcal{N}$ and we denote its basis with $\{ \phi_i \}^\mathcal{N}_{i=1}$. The fully \emph{discretized} problem reads
\begin{equation}\label{parabolic_problem_discretized}
\begin{split}
&\text{for each }1\leq j \leq J,\, \text{find }u^\mathcal{N}_j(\bmu(\omega)) \in X^\mathcal{N} \text{ s.t.}\\
&\frac{1}{\Delta t} m(u_j^\mathcal{N}(\bmu(\omega))-u_{j-1}^\mathcal{N}(\bmu(\omega)),v^\mathcal{N};\bmu(\omega))+a(u_j^\mathcal{N}(\bmu(\omega)),v^\mathcal{N};\bmu(\omega))=\\
&\qquad g(t_j)F(v^\mathcal{N};\bmu(\omega))\quad \forall v^\mathcal{N} \in X^\mathcal{N},\quad \forall \omega \in A,\\
&\text{given the initial condition } u^\mathcal{N}_0\text{ s.t.}\\
&(u^\mathcal{N}_0,v^\mathcal{N})_{L^2(\Omega)}=(u_0,v^\mathcal{N})_{ L^2(\Omega)}\quad \forall v^\mathcal{N}\in X^\mathcal{N}.
\end{split}
\end{equation}
The latter problem uses the \emph{Backward Euler-Galerkin} discretization, but we can resort to other theta-methods (e.g. Crank-Nicholson) or to high order method (e.g. Runge--Kutta) \cite{quarteroni_sacco}.

The RB formulation of the problem \eqref{parabolic_problem_discretized} is based on hierarchical RB space, as we did for the steady case, employing a POD reduction over the time trajectory and a greedy selection over the parameter space \cite{Haasdonk, Calcolo}. The algorithm can be seen as a Greedy algorithm in the parameter space with a further compression by POD for the space trajectory.

At each step of the Greedy algorithm we search the parameter $\bmu^*$ which maximizes, over the training set $\Xi_{train}$, an error estimator for the following quantity:
\begin{equation}\label{parabolic_error}
|||\bm{e}^\mathcal{N}_N(\bmu)|||_{t-dep}=\left( m(e_{N,J}^\mathcal{N}(\bmu), e_{N,J}^\mathcal{N}(\bmu);\bmu)+\sum_{j=1}^J a(e_{N,j}^\mathcal{N}(\bmu), e_{N,j}^\mathcal{N}(\bmu);\bmu)\Delta t \right)^\frac{1}{2},
\end{equation} 
where $e^\mathcal{N}_{N,j}(\bmu)=u^\mathcal{N}_j(\bmu)-u^\mathcal{N}_{N,j}(\bmu)$. We remark that, as in section \ref{RB}, an inexpensive \emph{a posteriori} error bound for \eqref{parabolic_error} can be derived (see \cite{grepl}), which in particular does not require any $\mathcal{N}$-dependent computation (e.g. it does not require the time trajectory to be computed for every $\bmu$ in the training set). We will continue denoting by $\Delta_N$ the resulting error estimator, even though its expression is different from the one in section \ref{RB}; we refer to \cite{grepl} for more details.

Once the parameter is chosen, we project the time evolution of the solution of this parameter on the orthogonal space of the current reduced basis space $X_N^\mathcal{N}$. This projection ensures that, at each Greedy iteration, only \emph{new} information is added to the reduced basis. To set the notation, denote by $\mathcal{P}_N:X^\mathcal{N}\to X_N^\mathcal{N}$ the projection onto the current reduced basis $X_N^\mathcal{N}$. We then define $u_j^\perp(\bmu^*) = u_j(\bmu^*)-\mathcal{P}_N(u_j(\bmu^*))$, for $j = 1, \hdots, J$.

As a further compression of the resulting time trajectory, we compute a POD on $\lbrace u^\perp_j (\bmu^*)\rbrace_{j=1}^J$, and collect the first few POD modes (up to a prescribed tolerance) into a space denoted by $Y_N^\mathcal{N}$. The resulting reduced basis space to be used at the ($N+1$)-th Greedy iteration is then defined as $X_{N + 1}^\mathcal{N} = X_N^\mathcal{N} \oplus Y_N^\mathcal{N}$. 

The RB formulation of the problem can be obtained by substituting the reduced basis space $X^\mathcal{N}_N$ to $X^\mathcal{N}$ in \eqref{parabolic_problem_discretized}.

\subsection{SUPG stabilization method for parabolic problems}
In this section we briefly introduce the SUPG method for time-dependent problems \cite{brooks, johnson2}. The idea is the same of the steady case: we have to add terms to bilinear forms in order to improve stability. The stabilization term is almost the same than in the steady case, but now we have to consider also the time dependency to guarantee the strong consistency.
We thus set
\begin{equation}\label{stabilization_parabolic_lhs}
s(w^\mathcal{N}(t),v^\mathcal{N})=\sum_{K\in \mathcal{T}_h} \delta_K \left( \partial_t w^\mathcal{N}(t) + Lw^\mathcal{N}(t), \frac{h_k}{|\bbe (\bmu(\omega))|}L_{SS}v^\mathcal{N}\right)_K
\end{equation}
where $w^\mathcal{N}(t)\in X^\mathcal{N}$ for each $t\in I$, $v^\mathcal{N}\in X^\mathcal{N}$ and $\left(\cdot,\cdot \right)_K$ is the usual $L^2$ scalar product, restricted to the element $K$. Here $L$ is the steady advection--diffusion operator and $L_{SS}$ is its skew--symmetric part.
\par
%We note that if either the coefficients of the equation or its domain are $\bmu$-dependent, then the stabilization terms will depend on $\bmu$ too.
Thus, we can  define the \emph{Backward Euler--SUPG} formulation of the problem by substituting the forms $m$, $a$ and $F$ in \eqref{parabolic_problem_discretized} with:
\begin{equation}\label{parabolic_stabilized_terms}
\begin{split}
&m_{stab}(w^\mathcal{N},v^\mathcal{N};\bmu(\omega))=m(w^\mathcal{N},v^\mathcal{N};\bmu(\omega)) + \sum_{K\in \mathcal{T}_{h}} \delta_{K} \left(w^\mathcal{N}, \frac{h_{K}}{|\bbe(\bmu(\omega))|}L_{SS} v^\mathcal{N} \right)_{K}\\
&a_{stab}(w^\mathcal{N},v^\mathcal{N};\bmu(\omega))=a(w^\mathcal{N},v^\mathcal{N};\bmu(\omega)) + \sum_{K\in \mathcal{T}_{h}} \delta_{K} \left( Lw^\mathcal{N}, \frac{h_{K}}{|\bbe(\bmu(\omega))|}L_{SS} v^\mathcal{N} \right)_{K}\\
&F_{stab}(v^\mathcal{N};\bmu(\omega))=F(v^\mathcal{N};\bmu(\omega)) + \sum_{K\in \mathcal{T}_{h}} \delta_{K} \left( f, \frac{h_{K}}{|\bbe(\bmu(\omega))|}L_{SS} v^\mathcal{N} \right)_{K}
\end{split}
\end{equation}
where $K$ are the elements which form the mesh $\mathcal{T}_{h}$ and $f$ can be a source term of the advection--diffusion equation or a lifting of the Dirichlet boundary data. For the analysis of stability and convergence of the method we refer to \cite{John_Novo}.

\subsection{Numerical tests for stochastic parabolic problems}
\begin{figure}
\centering
\subfigure[$t=0$ sec]{\includegraphics[trim=0cm 5cm 1cm 5cm,clip=true, scale=0.35]{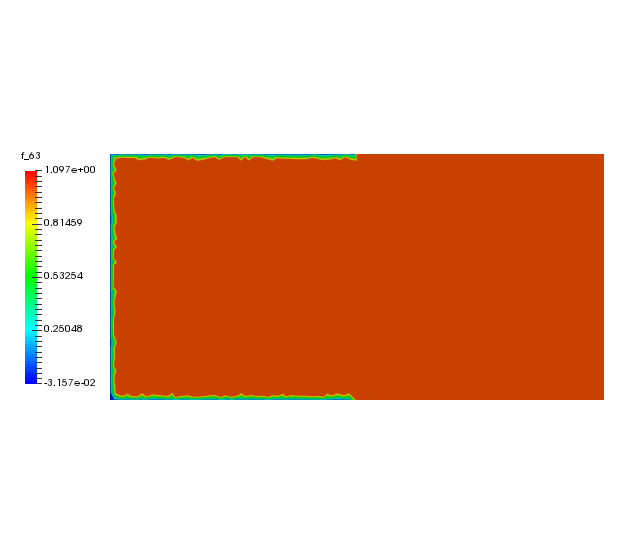} }
\subfigure[$t=1$ sec]{\includegraphics[trim=4cm 5cm 1cm 5cm,clip=true, scale=0.35]{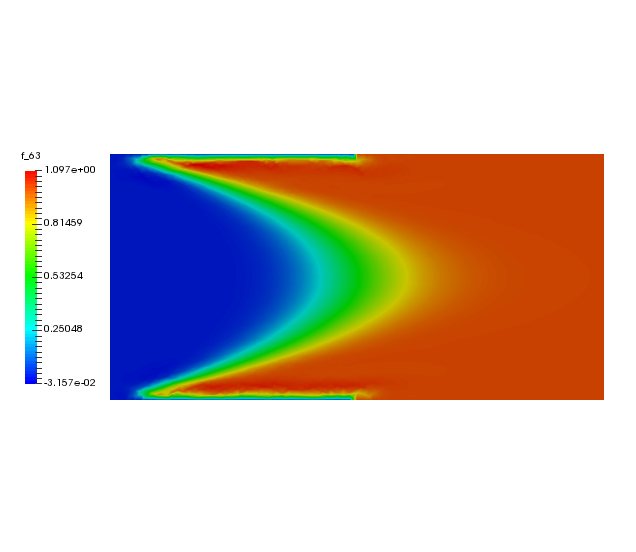} }\\
\subfigure[$t=2$ sec]{\includegraphics[trim=0cm 5cm 1cm 5cm,clip=true, scale=0.35]{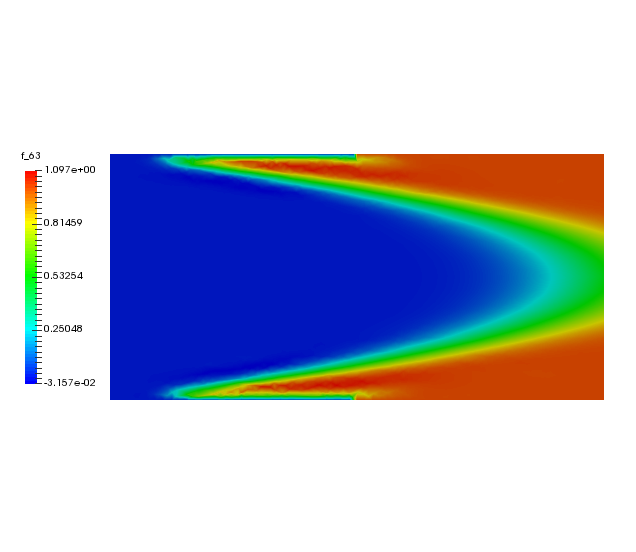} }
\subfigure[$t=7$ sec]{\includegraphics[trim=4cm 5cm 1cm 5cm,clip=true, scale=0.35]{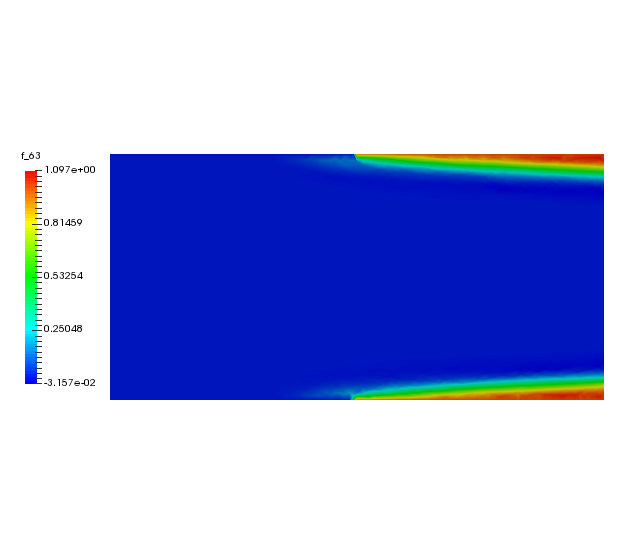} }\\
\caption{Plot of FE solution for parabolic \PGP problem at different times at $\mu_1=1$ and $\mu_2=1\cdot 10^4$}\label{parabolic_graetz_simulation}
\end{figure}

\begin{figure}
\centering
\subfigure[$t=0$ sec]{\includegraphics[trim=4cm 2cm 4cm 3cm,clip=true, scale=0.31]{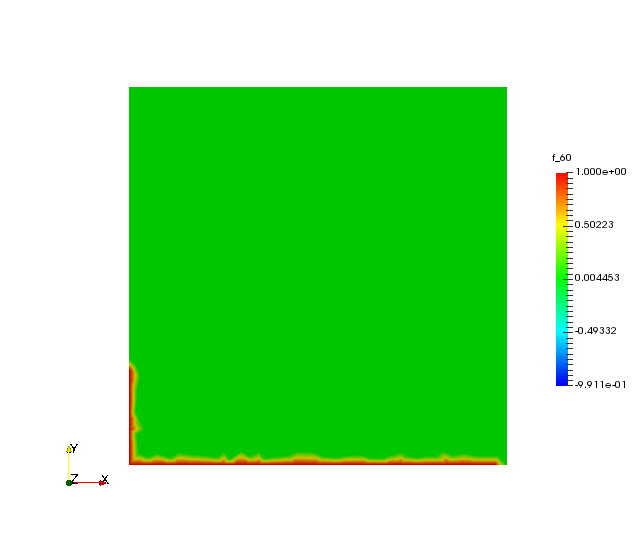} }
\subfigure[$t=0.64$ sec]{\includegraphics[trim=4cm 2cm 4cm 3cm,clip=true, scale=0.31]{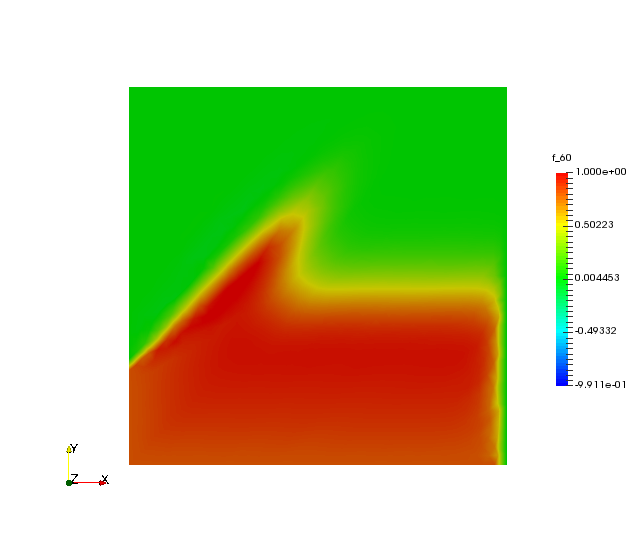} }
\subfigure[$t=1.28$ sec]{\includegraphics[trim=4cm 2cm 1cm 3cm,clip=true, scale=0.31]{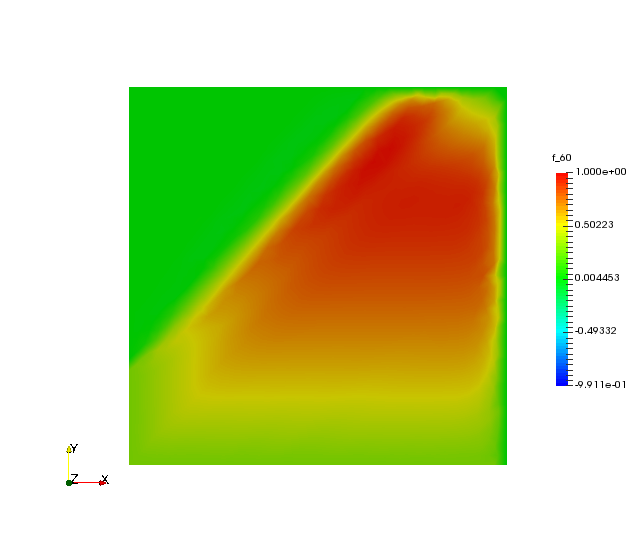} }
\caption{Plot of FE solution for parabolic \PFS problem at different times, $\mu_1=2\cdot 10^4,\,\mu_2=0.8$}\label{square_parabolic_simulation}
\end{figure}

We are now showing some numerical results of the stabilized RB method for stochastic parabolic PDEs, extending to the time dependent case the problems in sections \ref{stochastic_graetz} and \ref{stochastic_square}. For the sake of exposition we will show the results only for the Offline-Online stabilization. Few representative FE solutions are provided in figure \ref{parabolic_graetz_simulation} for the parabolic \PGP problem and figure \ref{square_parabolic_simulation} for the parabolic front propagation test.

\begin{figure}[h]
\centering
\subfigure[\PGP Problem]{\includegraphics[trim=0cm 0cm 1.5cm 0cm,clip=true, scale=0.38]{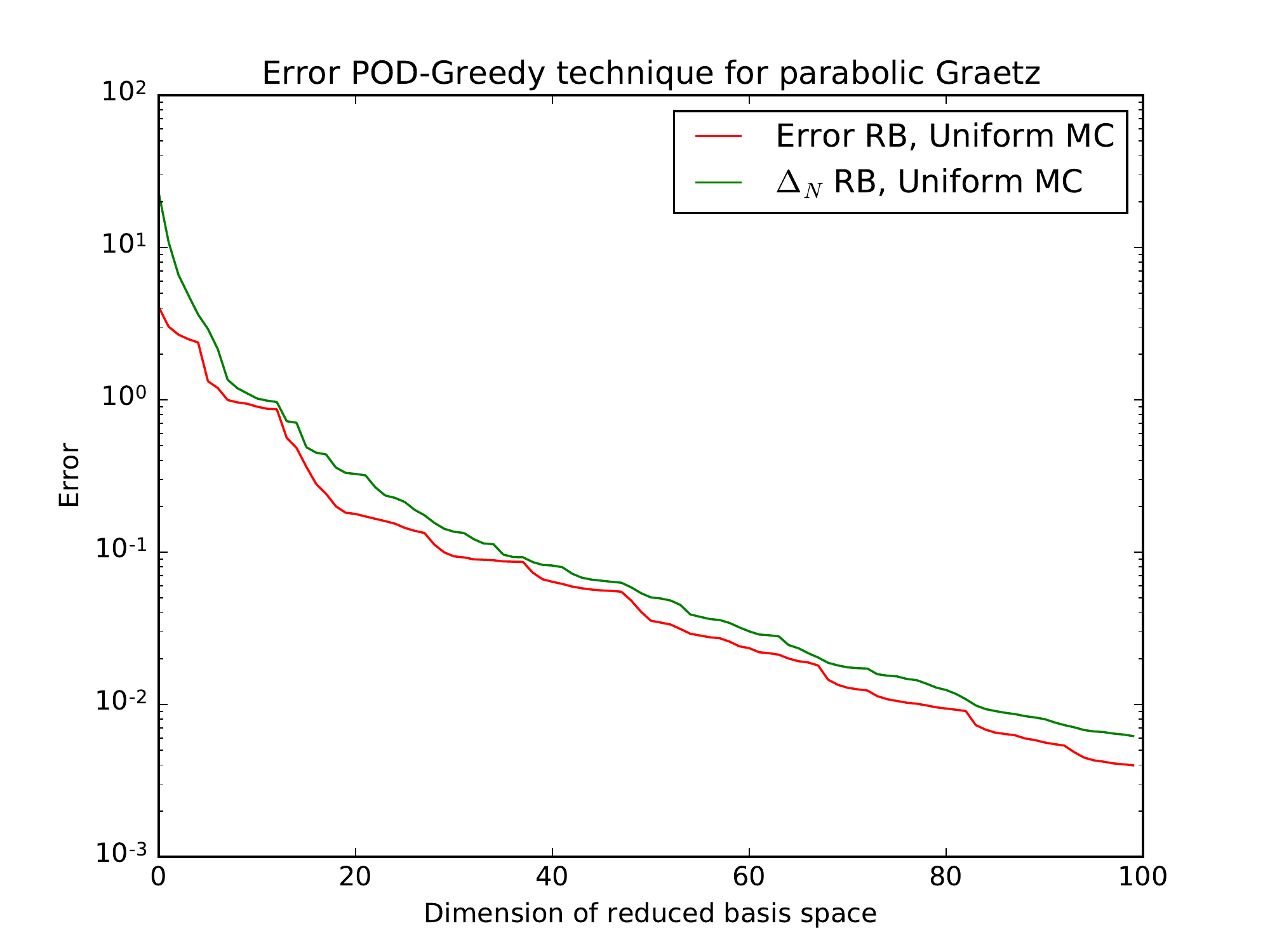} \label{graetz_parabolic}}
\subfigure[\PFS Problem]{\includegraphics[trim=0cm 0cm 1.5cm 0cm,clip=true, scale=0.38]{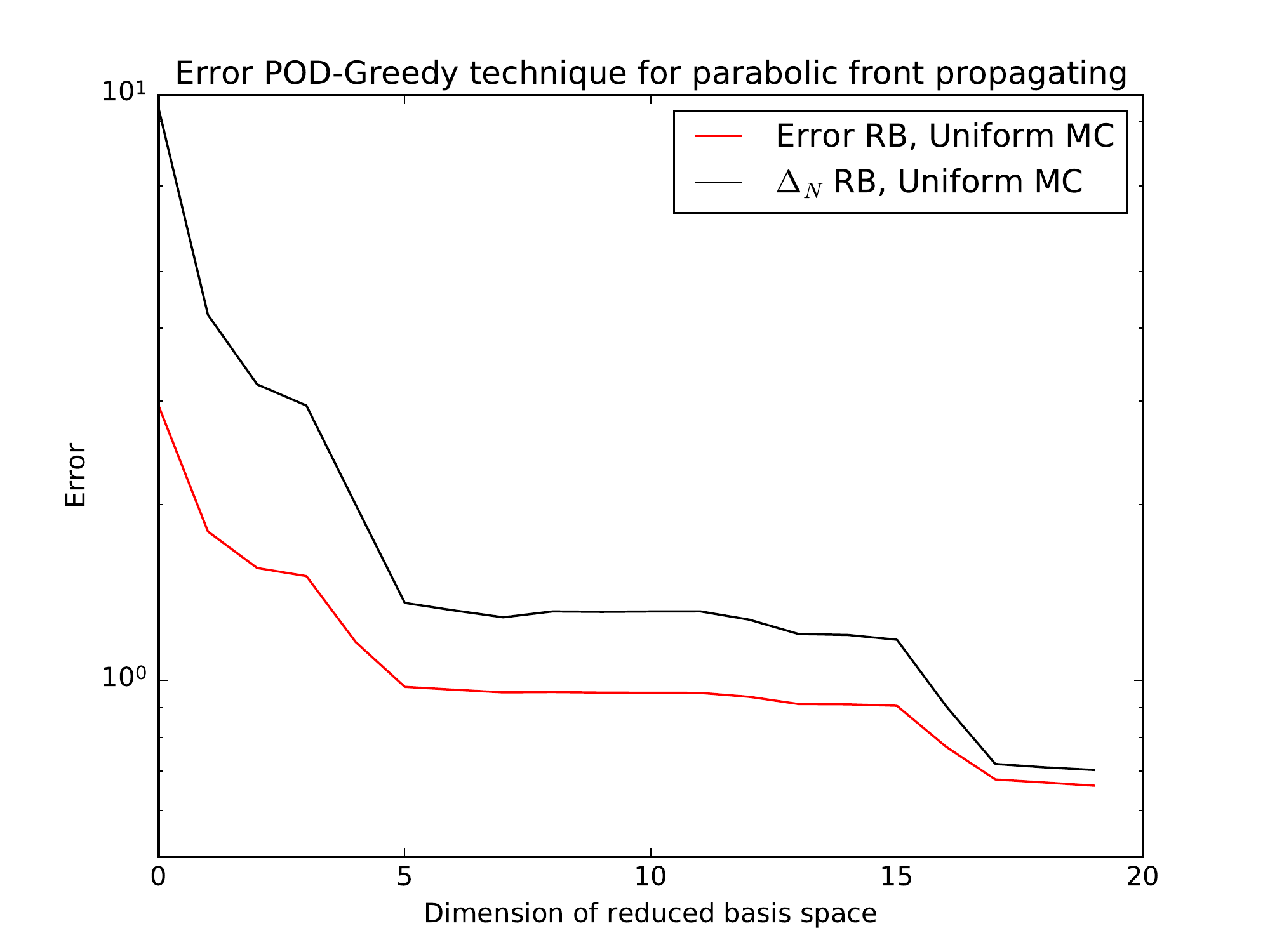} \label{square_parabolic} }
\caption{Greedy algorithms comparison for parabolic problems} \label{parabolic_greedy_errors}
\end{figure}

\begin{figure}[h]
\centering
\subfigure[\PGP Problem]{\includegraphics[trim=0cm 0cm 1.5cm 0cm,clip=true, scale=0.38]{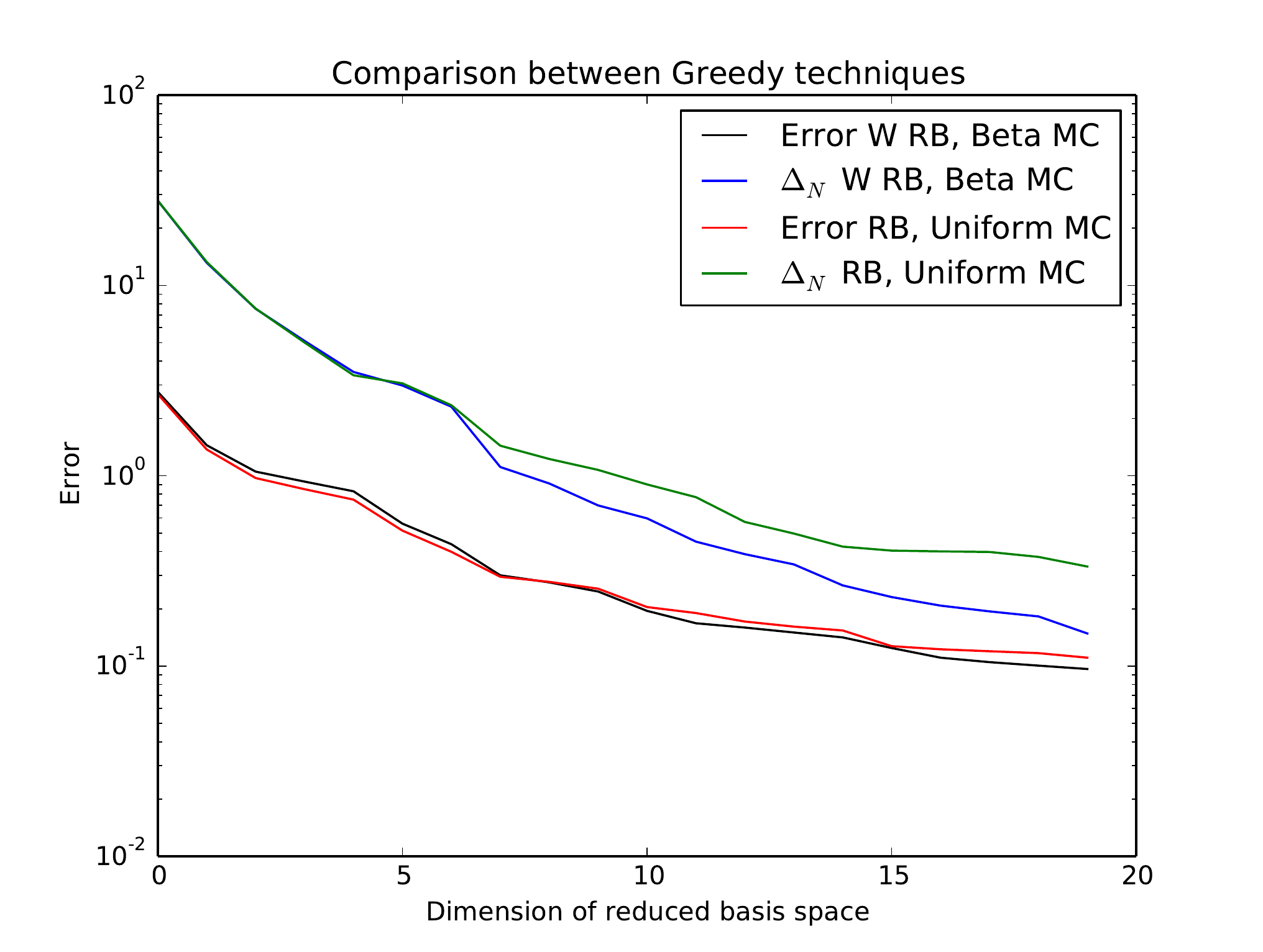} \label{stochastic_graetz_parabolic}}
\subfigure[\PFS Problem]{\includegraphics[trim=0cm 0cm 1.5cm 0cm,clip=true, scale=0.38]{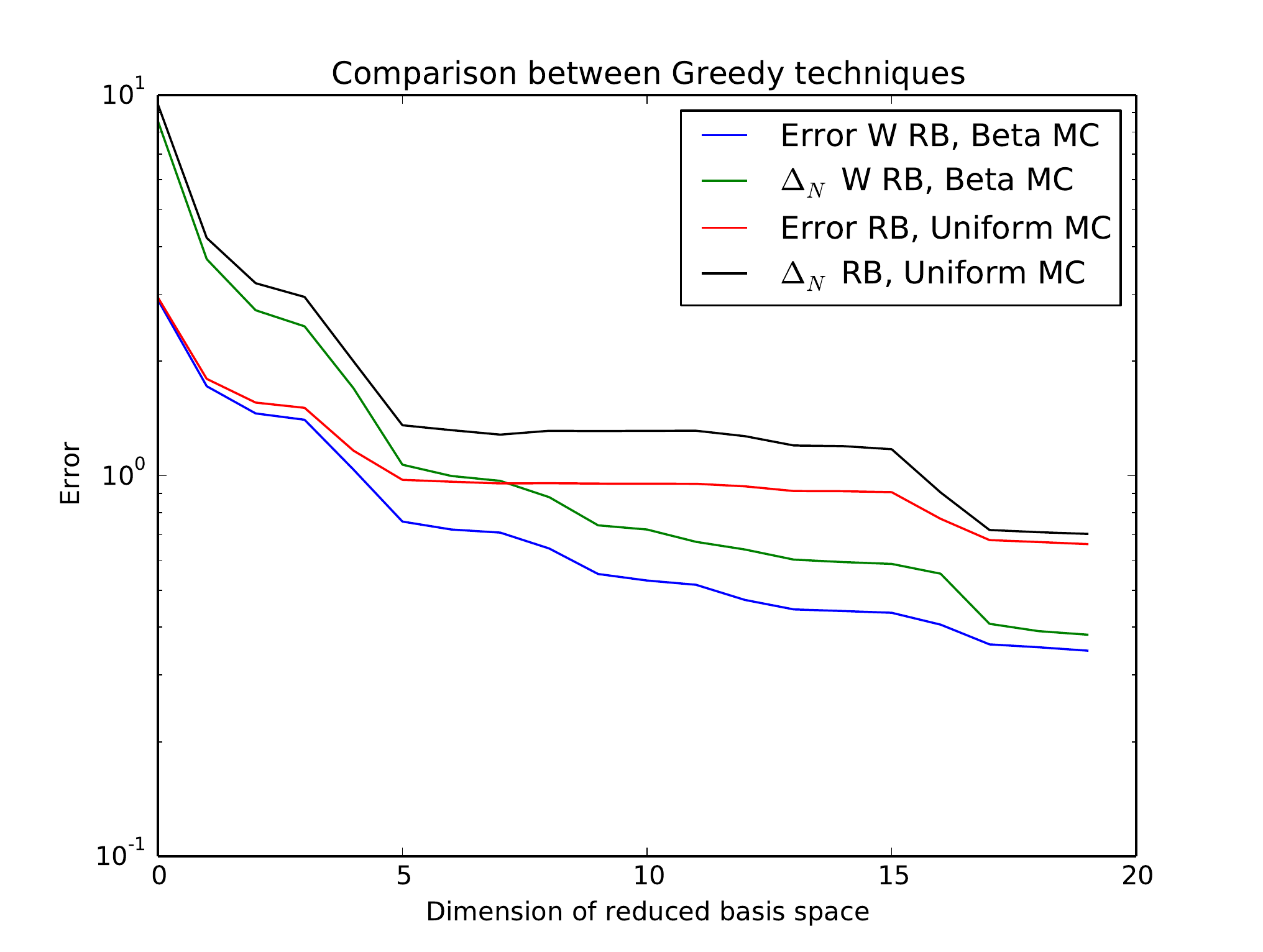} \label{stochastic_square_parabolic} }
\caption{Greedy algorithms comparison for parabolic problems} \label{stochastic_parabolic_comparison_greedy}
\end{figure}

We show in figures \ref{parabolic_greedy_errors} and \ref{stochastic_parabolic_comparison_greedy} the average error on a test set, for both the parabolic \PGP problem (left) and the parabolic front propagation test (right), respectively in deterministic and stochastic case. The error is the one defined in \eqref{parabolic_error}, while the error estimator $\Delta_N$ is as in \cite{grepl}. We compare in figure \ref{stochastic_parabolic_comparison_greedy} the classical reduced basis algorithm (with uniform Monte Carlo sampling) and the weighted reduced basis one (with sampling according to the distribution of $\bmu$). The comparison shows that, also for parabolic problem, proper weighting and suitable sampling allows to improve the accuracy of the resulting reduced order model (especially in the case of the parabolic front problem) and the reliability of the error estimator (in both test cases).

%We can see in figure \ref{stochastic_parabolic_comparison_greedy} that for Graetz problem the difference between classical Greedy method and weighted one is not so relevant for the first 20 dimension of the reduced basis space in terms of the error, while the $\Delta_N$ error bound is improved with the modified algorithm. For propagation front problem we can see important improvements both in error and $\Delta_N$.

Similar results hold for the probabilistic mean indicator introduced in \eqref{error_stochastic}, which we extend to the unsteady case as
\begin{equation}\label{error_parabolic_stochastic}
\mathbb{E}[||\bm{u}^\mathcal{N}-\bm{u}^\mathcal{N}_N||^2]:=\sum_{j=1}^J \int_\mathcal{D} ||u^\mathcal{N}_j(\bmu) -u_{N,j}^\mathcal{N}(\bmu)||^2_{\bmu} \rho(\bmu) d\bmu
\end{equation}
and approximate with Monte Carlo quadrature procedure. 
By doing this we obtain for \PGP problem with a reduced basis space of dimension 20 an error of $8.3248 \cdot 10^{-2}$ for classic Greedy algorithm and $7.6318 \cdot 10^{-2}$ for weighted reduced basis algorithm, respectively. For \PFS problem we have that the classic Greedy algorithm produce an error of $0.3196$ while the weighted algorithm gets $0.2343$. 
%\par In conclusion we can say that the weighted reduced basis method can be useful in situation where we know the distribution of parameters. Moreover, in the particular case where the parameters are concentrated around a high P\'eclet zone, we can use this algorithm to reduce our Online computation.
\par A small remark on computational times in parabolic must be done. In \PGP problem for one true parabolic solution we need 132.382 seconds, while for the RB one with $N=20$ basis functions we need only 0.356224 seconds. For a \PFS true solution we need 17.2846 seconds and only 0.125266 seconds for RB solution with $N=20$ basis functions. These results justify all the computational costs of the \textit{Offline} phase.  

\section{Conclusions}
\label{conclusions}
In this work we have dealt with stabilization techniques for the approximation of advection dominated problems using a reduced basis approach into a stochastic framework, both in steady and unsteady case. To perform a stabilization in the reduced basis algorithm, we have studied the SUPG \cite{quarteroni_valli} stabilization for FE method and introduced two reduced basis stabilization algorithms. The \emph{Online--Offline} stabilization, which uses SUPG stabilized forms in both stages (\emph{Offline} and \emph{Online}) and the \emph{Offline--only} stabilization, which uses the original (not stabilized) forms for the \emph{Online} stage. The underlying idea was to obtain a stable RB approximation, from the stable FE approximation, with reasonable computational times and, at the same time, a very good accuracy.
%\par We have tested the two methods on some examples which have shown that the \emph{Offline--Online} stabilization produces stable results and that the \emph{a posteriori} error bound of RB method still guarantees a convergence. The \emph{Offline--only} shows some instabilities and the error between RB approximation and FE one may not always converge because of declared inconsistencies between \emph{Online} and \emph{Offline} phase forms.
\par We then introduced stochastic equations and weighted reduced basis method \cite{peng_art}. We formulated a stabilized weighted reduced basis method for advection-diffusion problems with random input parameters. Numerical test cases clearly highlight the importance of the weighting procedure, as well as the necessity of a proper sampling of the parameter space, according to the probability distribution of $\bmu$.
Moreover, we introduced a procedure to selectively enable online stabilization when required. This allows to reduce the number of terms to be assembled in the affine expansion, with a negligible worsening of the error, which remains of the same order as the one for the previous strategies. 
\par Finally, we have generalized these methods to parabolic problems producing a stabilized RB approach for unsteady cases \cite{Haasdonk, pacciarini_a}, starting from SUPG stabilized parabolic FE methods \cite{brooks, johnson2}. 
%We have proposed two different algorithms to perform the parabolic RB algorithm and then we have compared them. We have tested successfully our algorithms on time dependent problems with also time dependent boundary conditions.

\par Possible further developments of this topic could be the application of these methods to more complex geometries, e.g. non--affinely parametrized ones, requiring some empirical interpolation preprocessing \cite{barrault, lassila_free_form}. Moreover, the method could be tested on larger dimension parameter spaces $\mathcal{D}$, using Monte Carlo or quasi--Monte Carlo strategies and on other types of probability distributions.

\section*{Acknowledgments} 
We acknowledge the support by European Union Funding for Research and Innovation -- Horizon 2020 Program -- in the framework of European Research Council Executive Agency: H2020 ERC Consolidator Grant 2015 AROMA-CFD project 681447 ``Advanced Reduced Order Methods with Applications in Computational Fluid Dynamics''. We also acknowledge the INDAM-GNCS projects ``Metodi numerici avanzati combinati con tecniche di riduzione computazionale per PDEs parametrizzate e applicazioni'' and ``Numerical methods for model order reduction of PDEs''. The computations in this work have been performed with RBniCS \cite{rbnics} library, developed at SISSA mathLab, which is an implementation in FEniCS \cite{fenics} of several reduced order modelling techniques; we acknowledge developers and contributors to both libraries.

%\nocite{*}
\bibliographystyle{siamplain}
\bibliography{biblio}

\end{document}